\documentclass[11pt]{article}
\usepackage{fullpage}
\usepackage{amssymb}
\usepackage{graphicx}

\begin{document}

\def\Title{Around Vogt's theorem}

\thispagestyle{empty}

\title{\Title}

\author{Alexey~Kurnosenko}

\maketitle{}


\centerline{\hbox{\it Institute for High Energy Physics, Protvino, 
     Russia,} 142281}

\begin{abstract}
Vogt's theorem, concerning boundary angles
of a convex arc with monotonic curvature ({\em spiral arc}), 
is taken as a starting point to establish basic properties
of spirals. The theorem is expanded by removing requirements
of convexity and curvature continuity; the cases of inflection 
and multiple windings are considered. 
Positional restrictions for a spiral arc with two given 
curvature elements at the endpoints are established,
as well as the necessary and sufficient conditions
for the existence of such spiral.
\\~~\\
\end{abstract}

\vfill
\noindent
{E-mails:            \verb"Alexey.Kurnosenko@ihep.ru"}\\
{\hphantom{E-mails:} \verb"Alexey.Kurnosenko@cern.ch"}
\thispagestyle{empty}

\newpage
\pagenumbering{arabic}

\newcommand{\tmp}{}
\newcommand{\Mark}[1]{${}^{#1}$\marginpar{{\small#1}}}

\newcommand{\Topic}[1]{\section{#1}}

\makeatletter
\def\@begintheorem#1#2{\trivlist \item[\hskip \labelsep{\bf #1\ #2.}]\it}
\def\@opargbegintheorem#1#2#3{\trivlist
      \item[\hskip \labelsep{\bf #1\ #2\ (#3).}]\it}
\makeatother
%
\newcommand{\Equa}[2]{\begin{equation}#2\label{#1}\end{equation}}
\newcommand{\equa}[1]{\[ #1 \]}

\newcommand{\refeq}[1]{{\rm(\ref{#1})}}
\newcommand{\refeqeq}[2]{{\rm(\ref{#1},\ref{#2})}}
\newcommand{\sign}{\mathop{\rm sign}\nolimits}
%
\newcommand{\iu}{{\mathrm{i}}}
\newcommand{\e}{{\mathrm{e}}}
%
\newcommand{\D}[1]{#1^{\prime}}
\newcommand{\DD}[1]{#1^{\prime\prime}}
\newcommand{\Int}[4]{\displaystyle\int\limits_{#1}^{#2}{#3}\,{\mathrm{d}}#4}
\newcommand{\Dfrac}[2]{\Frac{{\rm d}#1}{{\rm d}#2}}
\newcommand{\Pd}[2]{\Frac{\partial#1}{\partial#2}}
%
%
\newcommand{\ieq}{\,{=}\,}
\newcommand{\ineq}{\,{\not=}\,}
\newcommand{\ilt}{\,{<}\,}
\newcommand{\igt}{\,{>}\,}
\newcommand{\In}{\,{\in}\,}
\renewcommand{\le}{\leqslant}
\newcommand{\ile}{\,{\le}\,}
\renewcommand{\ge}{\geqslant}
\newcommand{\ige}{\,{\ge}\,}
\newcommand{\eqref}[1]{\stackrel{\refeq{#1}}{=}}
%
\newcommand{\kn}{\kappa}   

\newcommand{\Arc}[1]{\displaystyle{\buildrel\,\,\frown\over{#1}}}
\newcommand{\Biarc}{{\cal B}}
\newcommand{\Biarcab}[3]{\Biarc(#1;#2,#3)}
\newcommand{\Aarc}[1]{{\cal A}(#1)}
\newcommand{\Lense}[1]{{\mathbf L}(#1)}
\newcommand{\Bilense}[1]{{\bf B}(#1)}
%
%
\newcommand{\GT}[1]{$#1\igt0$}
\newcommand{\GE}[1]{$#1\ige0$}
\newcommand{\LT}[1]{$#1\ilt0$}
\newcommand{\LE}[1]{$#1\ile0$}
\newcommand{\EQ}[1]{$#1\ieq0$}
\newcommand{\NE}[1]{$#1\ineq0$}
%
\newcommand{\Frac}[2]{\displaystyle\frac{#1}{#2}}
\newcommand{\So}{\quad\Longrightarrow\quad}
\newcommand{\Vec}[1]{\stackrel{\longrightarrow}{{#1}}}
\newcommand{\HM}{\hphantom{{-}}}
\newcommand{\gr}[1]{#1\ifmmode^{\circ}\else$^{\circ}$\fi}

\newcommand{\Sa}{\sin\alpha}
\newcommand{\Sb}{\sin\beta}
\newcommand{\Mat}[1]{\mathop{\rm Mat}\nolimits(#1)}

\newcommand{\acum}{\widetilde\alpha}
\newcommand{\bcum}{\widetilde\beta}
\newcommand{\tcum}{\widetilde\tau}
\newcommand{\ocum}{\widetilde\omega}
\newcommand{\mcum}{\widetilde\mu}
\newcommand{\binf}{b^{\star}}

\newcommand{\Kl}[1]{{\ifmmode{\cal K}_{#1}\else${\cal K}_{#1}$\fi}}
\newcommand{\DKl}[1]{{\ifmmode{\cal K}^{\prime}_{#1}%
                      \else${\cal K}^{\prime}_{#1}$\fi}}
\newcommand{\Kr}[1]{K(#1)}
\newcommand{\Krn}[1]{\Kr{x_{#1},y_{#1},\tau_{#1},k_{#1}}}

%
%
\newcommand{\Figref}[1]{{\rm\ref{F#1}}}
\newcommand{\Reffigs}[1]{Figs.$\:$\Figref{#1}}
\newcommand{\Reffig}[1]{Fig.$\:$\Figref{#1}}
\newlength{\tmplength}                    %
\newcommand{\Infig}[3]{
\includegraphics[width=#1]{#2.eps}\caption{#3}\label{F#2}}
%
%
\newcommand{\Lfig}[2]{
\begin{figure}[#1]
\settowidth{\tmplength}{Fig.~9999.}%
\parbox[b]{\tmplength}{\caption{}\label{F#2}}%
\addtolength{\tmplength}{-\textwidth}%
\includegraphics[width=-\tmplength]{#2.eps}
\end{figure}
}
\newcommand{\Lwfig}[3]{
\begin{figure}[#1]%
\settowidth{\tmplength}{Fig.~9999.}%
\parbox[b]{\tmplength}{\caption{}\label{F#3}}%
\ifdim#2=0pt
  \addtolength{\tmplength}{-\textwidth}%
  \setlength{\tmplength}{-\tmplength}%
\else
  \setlength{\tmplength}{#2}%
\fi
\includegraphics[width=\tmplength]{#3.eps}
\end{figure}%
}
\newcommand{\Bfig}[3]{
\parbox[b]{#1}{\Infig{#1}{#2}{#3}}%
}

\newcommand{\Ffig}[3]{
   \begin{figure}[#1]
      \Infig{\textwidth}{#2}{#3}
   \end{figure}
}

 \newtheorem{thm}{Theorem}[section]

\newcommand{\qed}{\ifmmode{\qquad\mbox{\underline{q.e.d.}}}%
                  \else{{}~\hfill\underline{q.e.d.}\\}\fi} 
\newtheorem{lem}[thm]{Lemma}
\newtheorem{defn}[thm]{Definition}
\newtheorem{prop}[thm]{Proposition}
\newtheorem{cor}[thm]{Corollary}

\def\figurename{Fig.}

\makeatletter


\long\def\@makecaption#1#2{}
\newdimen\localfigtabsize
\localfigtabsize=0pt

\long\def\@makefigcaption#1#2{
 \ifdim \localfigtabsize > 0pt  \else \localfigtabsize=\hsize \fi
 \vskip 10pt
 \setbox\@tempboxa\hbox{\small#1.\quad#2}
 \ifdim \wd\@tempboxa > \localfigtabsize
 \settowidth{\@tempdima}{\small#1.\quad}
 \addtolength{\@tempdima}{-\localfigtabsize}
 {\small #1.\quad\parbox[t]{-\@tempdima}{#2}\par~\par}
 \else
 \hbox to\hsize{\hfil\box\@tempboxa\hfil}
 \fi
}

\let\@makecaption\@makefigcaption
\makeatother

\makeatletter
\newdimen\@bls
\@bls=.7\baselineskip      
\newenvironment{pf}%
  {\par\addvspace{\@bls plus 0.5\@bls
    minus 0.1\@bls}\noindent
   {\bf\proofname}\enspace\ignorespaces}%
  {\par\addvspace{\@bls plus 0.5\@bls minus 0.1\@bls}}
\def\proofname{Proof.}

\def\section{\@startsection{section}{1}{\z@}{1.5\@bls
  plus .4\@bls minus .1\@bls}{\@bls}{\normalsize\bf}}

\makeatother


\Topic{Introduction}
%
Vogt's theorem was published in 1914 (\cite{Vogt}, Satz~12).
It concerns convex arcs of planar curves with continuous
monotonic curvature of constant sign.
The later proofs~\cite{Tohoku1,Tohoku2,Ostrowski}
did not extend the class of invoked curves.
Guggenheimer (\cite{Guggen}, p.{\,}48) applies the term
{\em spiral arc} to such curves,
and formulates Vogt's theorem as follows:

{\em ``Let $A$ and $B$ be the endpoints of a 
spiral arc, the curvature nondecreasing from $A$ to $B$.
The angle~$\beta$
of the tangent to the arc at $B$ with the chord~$AB$
is not less than the angle~$\alpha$ of the tangent at~$A$ with~$AB$.
$\alpha\ieq\beta$ only if the curvature is constant''.}

Below $\alpha$ and $\beta$ denote algebraic values of the boundary angles 
with respect to the positive direction of $X$-axis, 
same as the direction of the chord~$\Vec{AB}$.
Signed curvatures at the endpoints~$A$ and~$B$ 
are denoted as~$k_1$ and~$k_2$.
Vogt assumes positive values for angles and curvatures,
and the theorem states that
$|\alpha|\igt|\beta|$ for $|k_1|\igt |k_2|$, and vice versa.
Five cases, depicted in \Reffig{figure01}
as arcs $AB_{i}$, can be detailed as
\equa{%
\begin{array}{lllcl}
{AB_1:}\quad&
\hphantom{0 > {}}k_1 < k_2 < 0,\quad&
{ \HM{}\alpha = |\alpha| > |\beta|= {-}\beta} &
\So&  \alpha{+}\beta > 0;\\[2pt]
{AB_2:}&
0 > k_1 >  k_2, &
{ \HM{}\alpha = |\alpha| < |\beta| = {-}\beta} &
\So&  \alpha{+}\beta  <  0;\\[2pt]
AB_0: &
\hphantom{0 > {}}k_1=k_2,&
{\pm}\alpha = |\alpha| = |\beta| = {\mp}\beta &
\So& \alpha{+}\beta = 0;\\[2pt]
{AB_3:}&
\hphantom{0 > {}} k_1 >  k_2 > 0, &
{ -\alpha =|\alpha| > |\beta| = \HM{}\beta} &
\So&  \alpha{+}\beta  <  0;\\[2pt]
{AB_4:}&
0 < k_1 < k_2,&
{ {-}\alpha =|\alpha|<|\beta| = \HM{}\beta} &
\So&  \alpha{+}\beta  >  0;
\end{array}
}
and unified to
\Equa{VogtAK}{%
     \sign(\alpha{+}\beta) = \sign(k_2{-}k_1)
}
for whichever kind of monotonicity and curvature sign.
In this notation the theorem remains valid for non-convex arcs.
The proof for this case, lemma~1 in~\cite{Spiral},
required curve to be one-to-one projectable onto its chord.

A variety of situations is illustrated by a family
of arcs of Cornu spirals with fixed $\alpha$
and varying $\beta$, shown in \Reffig{figure02}.
The sum $\alpha{+}\beta$ for curves~1 and~2 with decreasing curvature
is negative; it vanishes at circular arc~3, and becomes positive
for curves \hbox{4--10} with increasing curvature.
Vogt's theorem covers cases \hbox{1--4} (arc~5 is not convex), 
lemma~1 in~\cite{Spiral}~--- cases \hbox{3--8},
curves \hbox{3--4} are covered by both,
and \hbox{9--10}~--- by neither.
General proof for cases \hbox{1--10}, defined below as
{\em short} arcs, is the first extension
of Vogt's theorem proposed herein.

The requirements of convexity or projectability both
served to somehow shorten the arc.
But it turned out that Vogt's theorem can also be formulated for
``long'' spirals like curve~11 in \Reffig{figure02}.


\Topic{Preliminary definitions and notation} 
%
We describe curves by intrinsic equation $k\ieq k(s)$,
$k$~being curvature, and $s$, arc length:
$0\ile s\ile S$. Functions $Z(s)\ieq x(s){+}\iu y(s)$ and
$\tau(s)$ represent coordinates and the direction of tangent.
The terms {\em ``increasing'' }  and {\em ``decreasing'',}
applied to any function $f(s)$,
are accompanied in this article by the adverb {\em ``strictly''{} } when
necessary; otherwise non-strict monotonicity  with
$f(s){\not\equiv}const$ is assumed.
\begin{defn}
\label{DefSpiral0}
\rm
{\em Spiral} is a planar curve of monotonic,
piecewise continuous curvature,
not containing the circumference of a circle.
Inflection and infinite curvatures at the endpoints
are admitted.
\end{defn}
\begin{figure}[t]
\Bfig{.4\textwidth}{figure01}{}%
\hfill%
\Bfig{.4\textwidth}{figure02}{}%
\end{figure}
\begin{defn}
\label{DefBiarc}
\rm
{\em Biarc} is a spiral, composed of two arcs of constant curvature 
(like arcs $AT_3B_3$ and $AT_4B_4$ in \Reffig{figure01}).
\end{defn}
\begin{defn}
\label{DefShort1}
\rm
An arc $\Arc{AB}$ is {\em short}, 
if its tangent never achieves
the direction $\displaystyle{\buildrel\,\,\to\over{BA}}$,
opposite to the direction of its chord,
except, possibly, at the endpoints.
\end{defn}
\begin{defn}
\label{DefShort2}
\rm
An arc $\Arc{AB}$ is {\em short}, 
if it does not intersect the complement of its chord
to the infinite straight line (possibly intersecting the chord itself).
\end{defn}
Def.~\ref{DefShort1} will be used until the equivalence
of definitions~\ref{DefShort1} and~\ref{DefShort2} is proven 
(corollary~\ref{EquivDefCor}).
The term {\em ``very short'' {}} will be sometimes used
to denote an arc, one-to-one projectable onto its chord 
(curves 3--8 in \Reffig{figure02}).

Guggenheimer uses terms {\em line element} to denote
a pair $(P,{\bf t})$ of a point and a direction,
and {\em curvature element} $(P,{\bf t},\rho{\bf n})$,
${\bf n}\perp{\bf t}$, adding curvature radius~$\rho$
at~$P$ (\cite{Guggen}, p.~50).
We modify these definitions to $(x,y,\tau)$ and $(x,y,\tau,k)$
with tangent angle~$\tau$ and signed curvature~$k$ 
at the point $P\ieq(x,y)$. Notation 
\equa{%
  \Kl{i}=\Krn{i}
}
serves to denote both the $i$-th curvature element
and {\em directed curve of constant curvature}
produced by \Kl{i}. Whether it be a straight line or a circular arc,
it goes under general name {\em circle} [{\em of curvature}].

As in~\cite{InvInv}, we use an implicit equation of the circle 
$\Kl{0}\ieq\Krn{0}$
in the form
\Equa{Cxy}{%
C(x,y;\Kl{0}) \equiv   
k_0\left[(x{-}x_0)^2{+}(y{-}y_0)^2\right]+
2(x{-}x_0)\sin\tau_0 - 2(y{-}y_0)\cos\tau_0=0.
}
The sign of $C(x,y;\Kl{0})$ reflects the position of
the point $(x,y)$ with respect to~\Kl{0}:
it is from the left (\LT{C}) or from the right (\GT{C})
of the circle's boundary. Keeping in mind applications to
geometric modelling, define the following:
\begin{defn}
\rm
The {\em region of material} of the circle $\Kl{}$ is
\equa{
       \Mat{\Kl{}}=\{(x,y):\ C(x,y;\Kl{})\le 0\}.
}
\end{defn}
\Lwfig{t}{0pt}{figure03}%
To consider properties of a spiral arc in relation to its chord 
$\Vec{AB}$, $|AB|\ieq 2c$,
we choose the coordinate system such that the chord
becomes the segment $[-c,c]$ of X-axis. With 
$\alpha\ieq\tau(0)$, $k_1\ieq k(0)$, and 
$\beta\ieq\tau(S)$, $k_2\ieq k(S)$,
the boundary circles of curvature take form
\Equa{K1K2cc}{%
  \Kl{1} = K(-c,0,\alpha,k_1),\qquad
  \Kl{2} = K(c,0,\beta,k_2).
}
It is often convenient to assume homothety
with the coefficient $c^{-1}$, and to operate
on the segment $[-1,1]$. 
The coordinates $x,y$ and curvatures $k$ become
normalized dimensionless quantities, corresponding to \ $x/c$, \ $y/c$ 
and \ $kc\ieq\kn$.
With such homothety applied, 
boundary circles~\refeq{K1K2cc} appear as
\Equa{K1K2c1}{%
  \Kl{1} = K(-1,0,\alpha,\kn_1),\qquad
  \Kl{2} = K(1,0,\beta,\kn_2).
}
\begin{defn}
\rm
A curve whose start point is moved into position $A(-c,0)$,
and the endpoint, into $B(c,0)$, is named {\em normalized arc}.
The product $ck(s)\equiv \kn(s)$, invariant under homotheties,
will be referred to as {\em normalized curvature}.
\end{defn}
Denote $\Aarc{\xi}$ a normalized circular arc,
traced from the point $A(-c,0)$ to $B(c,0)$
with the direction of tangents $\xi$ at~$A$
\ ($k_{1,2} \ieq {-}\sin\xi/c$, \ $\kn_{1,2} \ieq {-}\sin\xi$).
The arc $\Aarc{\pm\pi}$, passing through infinity, is coincident
with the chord's complement to an infinite straight line;
the arc $\Aarc{0}$ is the chord $AB$ itself.
\begin{defn}
\rm A {\em lense} $\Lense{\xi_1,\xi_2}$ is the region
between
two arcs $\Aarc{\xi_1}$ and $\Aarc{\xi_2}$, namely
\equa{%
 \Lense{\xi_1,\xi_2}= \{\,(x,y):\; (x,y)\in\Aarc{\xi}\,\},\quad
   \min(\xi_1,\xi_2) < \xi < \max(\xi_1,\xi_2).
}
\end{defn}
The arc $\Aarc{\alpha}$ shares tangent with the
normalized spiral at the start point;
so does $\Aarc{-\beta}$ at the endpoint. 
The two arcs bound the lense
$\Lense{\alpha,-\beta}$, shown in gray in \Reffig{figure03}.
The signed half-width $\omega$ of the lense, and the direction~$\gamma$ 
of its {\em bisector} $\Aarc{\gamma}$
are
\Equa{OmegaGamma}{%
   \omega=\Frac{\alpha{+}\beta}{2},\qquad
   \gamma=\Frac{\alpha{-}\beta}{2}\,.
}
\begin{defn}
\rm
By the {\em inflection point} of a spiral, whose curvature $k(s)$
changes sign, shall be meant any inner point $Z(s_0)$ with \EQ{k(s_0)}.
If there is no such point, i.e. curvature jump
$k(s_0{-}0)\lessgtr 0 \lessgtr k(s_0{+}0)$ occurs,
the jump point will be used as the inflection point
with the assignment \EQ{k(s_0)}.
\end{defn}

\Topic{Vogt's theorem for short spirals}
%
The subsequent proof of the modified 
Vogt's theorem~\refeq{VogtAK} for short spirals is similar to the proof for ``very
short'' ones from~\cite{Spiral}. Both clearly show that only monotonicity of $k(s)$,
and not convexity of the arc, is the basis for Vogt's theorem.

\begin{thm}
\label{ShortVogtTheorem}
Boundary angles $\alpha$ and $\beta$
of a normalized short spiral or circular arc
obey the following conditions:
\Equa{VogtShort}{%
\begin{array}{llcc}
\mbox{if~~}k_1 <  k_2:\qquad& \alpha{+}\beta>0,\quad &
             -\pi < \alpha \le \pi,\; & -\pi < \beta  \le \pi; \\
\mbox{if~~}k_1 >  k_2:\qquad& \alpha{+}\beta<0,\quad & 
             -\pi \le \alpha < \pi,\; & -\pi \le \beta  < \pi;\\
\mbox{if~~}k_1 =  k_2:\qquad& \alpha{+}\beta=0,\quad &
             -\pi < \alpha < \pi,\; & -\pi < \beta  < \pi.
\end{array}
}
\end{thm}
\begin{pf}
Consider the case of increasing curvature $k(s)$,
and define a new parameter~$\xi$:
\equa{%
    \xi(s)=\Int{0}{s}{\cos\Frac{\tau(\sigma)}{2}}{\sigma},\quad
    0\le\xi\le\xi_1=\xi(S).
}
By Def.~\ref{DefShort1}, $|\tau(s)|\ilt\pi$ within
the interval $(0,\,S)$, and $\xi(s)$ is therefore strictly
increasing with~$s$. Define function  $z(\xi)$:
\equa{%
  z(\xi)=\sin\Frac{\tau(s(\xi))}{2}\,,\qquad
  \Dfrac{ z}{\xi}= \Dfrac{z}{s} \cdot \Dfrac{s}{\xi}=
  \left( \Frac{1}{2} \cos\Frac{\tau(s)}{2} \Dfrac{\tau}{s}\right)
  \cdot {\left(\cos\Frac{\tau(s)}{2}\right)}^{-1}
  =\Frac{1}{2}k(s(\xi)).
}
Its  derivative being increasing,
$z(\xi)$ is downwards convex, and its plot lies below straight line
segment $l(\xi)$, connecting the endpoints $z(0)\ieq\sin(\alpha/2)$
and $z(\xi_1)\ieq\sin(\beta/2)$:
\equa{%
  z(\xi)< l(\xi)=
  \Frac{\xi_1{-}\xi}{\xi_1}\sin(\alpha/2)+ \Frac{\xi}{\xi_1}\sin(\beta/2)\,,
  \qquad 0<\xi<\xi_1.
}
Condition $y(0)\ieq y(S)$ yields
\vspace{-10pt}
\equa{%
\begin{array}{rcl}
  0& = &\Int{0}{S}{\sin\tau(s)}{s}=
    2\int\limits_{0}^{S}\sin\Frac{\tau(s)}{2}
        \overbrace{\left(\cos \Frac{\tau(s)}{2}{\mathrm{d}}s\right)}^{d\xi}=
  2\Int{0}{\xi_1}{z(\xi)}{\xi} < \\
  & < &
  2\Int{0}{\xi_1}{l(\xi)}{\xi} = \Frac{\xi_1}{2}(\sin(\alpha/2)+\sin(\beta/2))
  =\xi_1\sin\Frac{\alpha{+}\beta}{4}\cos\Frac{\alpha{-}\beta}{4}\,.
\end{array}
}
So, inequality $|{\alpha{+}\beta}|\ile 2\pi$,
resulting from $\alpha,\beta\In[-\pi,\pi]$,
can be refined to 
\ $0\ilt\alpha{+}\beta\ile 2\pi$.
Condition \GT{\alpha{+}\beta}
excludes the value $-\pi$ for $\alpha$ or $\beta$, 
providing inequalities~\refeq{VogtShort} for the case
of increasing curvature.

If the curvature of the arc $Z(s)$ is decreasing,
i.e. $k_2\ilt k_1$, consider the arc $\bar{Z}(s)$,
symmetric to $Z(s)$ about $X$-axis.
Its boundary angles are
$\D{\alpha}\ieq{-\alpha}$, \ $\D{\beta}\ieq{-\beta}$,
and the curvature increases: $\D{k_1}\ieq{-}k_1 < {-}k_2=\D{k_2}$. So,
\GT{\D{\alpha}{+}\D{\beta}}, and
\LT{\alpha{+}\beta} for the original curve.
The case of constant curvature is evident.\qed{}
\end{pf}

\Topic{Basic inequality of the theory of spirals}
%
In article~\cite{InvInv} we have introduced
{\em an inversive invariant of a pair of circles}
$\Kl{1}\ieq\Krn{1}$ and $\Kl{2}\ieq\Krn{2}$:
\Equa{DefQ}{%
\begin{array}{rcl}
Q(\Kl{1},\Kl{2})
  &=&\Frac{1}{4}{ k_1 k_2}[(\Delta x)^2{+}(\Delta y)^2]+
                         \sin^2\Frac{\Delta\tau}{2}+{}\\[6pt]
  &+&\Frac{1}{2}{ k_2}(\Delta x\sin\tau_1{-}\Delta y\cos\tau_1)
     -\Frac{1}{2}{ k_1}(\Delta x\sin\tau_2{-}\Delta y\cos\tau_2)\\[6pt]
   &\HM{}& (\Delta x = x_2{-}x_1,\quad
       \Delta y = y_2{-}y_1,\quad \Delta \tau = \tau_2{-}\tau_1).
\end{array}
}
Its value is independent of arbitrarily chosen line elements
($x_i,y_i,\tau_i$) on each circle, and is invariant
under motions, homothety and inversions;
and $Q(\Kl{1},\Kl{2})=Q(\Kl{2},\Kl{1})$.
In particular cases,
\begin{enumerate}
\item[a)]
\NE{k_{1,2}}, $D$ is the distance between two centres;
\item[b)] 
\EQ{k_1}, \NE{k_2}, $L$ is (signed) distance from the centre 
of~\Kl{2} to the straight line~\Kl{1}; 
\item[c)]
any $k_1, k_2$, and $\psi$ is intersection angle of two circles;
\end{enumerate}
the invariant $Q$ can be represented as follows:
\Equa{Qabc}{%
Q^{(a)}=
   \Frac{(k_1k_2D)^2 - (k_2{-}k_1)^2}{4k_1k_2},\qquad 
Q^{(b)}=
   \Frac{1{+}k_2L}{2},\qquad 
Q^{(c)}=
   \sin^2(\psi/2).
}
If two circles have no real common point, $\psi$ is complex,
but $Q$ remains real. In this case\linebreak
\hbox{${\mathrm{Im}}(\psi)\ieq {\mathrm{cosh}}^{-1}|1{-}2Q|$} is Coxeter's inversive
distance of two circles~\cite{InvDist}.
\EQ{Q} if and only if two circles are tangent,
or are two equally directed straight lines.
Situation $Q\ieq 1$ can be considered as ``antitangency''
($\tau_1\ieq\tau_2\pm\pi$ at the common point).
Using this invariant is the alternative
to cumbersome enumeration of variants with different
mutual position (and curvature sign) of the two circles,
commonly occurring in describing problems of this sort.
%
%
\begin{prop}
\label{PropDeviate}
A curve with
increasing curvature intersects its every circle of curvature
from right to left (and from left to right if curvature decreases).
\end{prop}
Proving here this familiar statement helps us to have 
the subsequent proof  of our basic theorem~\ref{MainTheorem}
self-contained.  The second needed assumption,
invariance of ~$Q$,
is easy to prove by deducing formulae~\refeq{Qabc}.
\begin{pf}
To consider behavior of a spiral at some point $P_1\ieq Z(s_1)$
choose the coordinate system with the origin at $P_1$
and the axis of~X aligned with $\tau(s_1)$.
The curvature element at $P_1$ becomes then $\Kl{1}\ieq\Kr{0,0,0,k_1}$.
Obtain function $C(s)$ by substituting $x\ieq x(s)$ and $y\ieq y(s)$
into implicit equation~\refeq{Cxy} of~\Kl{1}:
\Equa{Cs}{%
    C(s)= C(x(s),y(s);\,\Kl{1})=k_1[x(s)^2+y(s)^2]-2y(s).
}
Differentiating $C(s)$ yields 
($x$, $y$, $\tau$ and $k$ are abbreviated functions of $s$):
\equa{%
  \begin{array}{lcl}
     C^{\prime}(s) &=& 
     2k_1(x\cos\tau+y\sin\tau)-2\sin\tau,\\
     C^{\prime\prime}(s) &=& 
        2k_1 + 
        2k_1 k(y\cos\tau - x\sin\tau) -
        2k\cos\tau,\\
     C^{\prime\prime\prime}(s) &=& 
        2k^2[\sin\tau-k_1(x\cos\tau{+}y\sin\tau)] -
        2\D{k}(s)[\cos\tau+k_1(x\sin\tau{-}y\cos\tau)];\\
     C(s_1) &=& 0,\qquad 
     C^{\prime}(s_1) = 0,\qquad
     C^{\prime\prime}(s_1) = 0,\qquad
     C^{\prime\prime\prime}(s_1) = {-}2k^\prime(s_1)\lessgtr 0.
  \end{array}
}
Whether $k(s)$ varies continuously or with a jump at $s_1$,
function $C(s)$ undergoes variation of the opposite sign;
$C(s)$ going negative (or positive, if curvature decreases), 
the curve locally deviates to the left (or to the right) 
of~\Kl{1}.\,\qed
\end{pf}
\begin{thm}
\label{MainTheorem}
Let \Kl{1} and \Kl{2} be two circles of curvature of a spiral curve.
Then
\Equa{Qle0}{%
             Q(\Kl{1},\Kl{2})\le 0,
}
and equality holds if and only if both circles
belong to a circular subarc or to a biarc.
\end{thm}
\begin{figure}[bt]
\Bfig{.3667\textwidth}{figure04}{}      
\hfill
\Bfig{.5333\textwidth}{figure05}{}     
\end{figure}
\begin{pf}
Denote the circle of curvature at the start point as \Kl{1}, 
$\Kl{}(s)$ being any other circle of curvature:
\equa{
  \Kl{}(s)= K(x(s),y(s),\tau(s),k(s)\,),\qquad
  \Kl{1}= \Kl{}(0).
}
Two examples in \Reffig{figure04} illustrate the proof for the case of
increasing curvature with negative and positive start values $k_1\ieq k(0)$.
The regions $\Mat{\Kl{1}}$ are shown in gray.
Points $P_i$ subdivide the spiral into subarcs
$0\ile s_1\ile s_2 \ile s_3 \ile s_4 \ile S$, some of them possibly absent:
$P_0P_1$ is the initial subarc of constant curvature (if any),
coincident with \Kl{1}.
As soon as the curvature increases at $P_1$, with or without jump,
the spiral deviates to the left from the circle~\Kl{1} (Prop.~\ref{PropDeviate}).
The arc $P_0P_1$ may be supplemented to a biarc by another
circular arc $P_1P_2$; point $P_3$ represents any curvature jump,
where there are two circles of curvature $\Kl{}(s_3{\pm}0)$,
left and right. Point $P_4\ieq Z(s_4)$, if exists, 
is the point, where local property~\ref{PropDeviate} is no more valid, 
i.e.  spiral returns to the boundary of~\Kl{1}. Thus, with expression~\refeq{Cs} 
for $C(s)$ involved,
\equa{%
  \begin{array}{l}
    C(s) = 0\mbox{~~~~if~~~} 0 \le s\le s_1, \mbox{~~~or~~} s=s_4,\\
    C(s) < 0\mbox{~~~~if~~~} s_1 < s < s_4.                    
  \end{array}
}

Associate the coordinate system with the line element
$(x(s_1),y(s_1),\tau(s_1)\,)$ such that\linebreak
$\Kl{1}\ieq\Kl{}(0)\ieq\Kl{}(s_1{-}0)\ieq\Kr{0,0,0,k_1}$.
Define function $Q(s)= Q(\Kl{1},\Kl{}(s))$ from~\refeq{DefQ}:
\Equa{Qs}{%
    Q(s)=
    \Frac{1}{4}k(s)C(s)
    -\Frac{k_1}{2}[x(s)\sin\tau(s){-}y(s)\cos\tau(s)]
    +\sin^2\Frac{\tau(s)}{2}\qquad [Q(0)=0]
}
Show that $Q(s)$ is monotonic decreasing in $[0,s_4]$. 
Differentiating of $Q(s)$ yields
\Equa{dQs}{%
  \D{Q}(s) = \Frac{1}{4}\D{k}(s)\,C(s)\So \D{Q}(s)  \le 0
}
for increasing $k(s)$. Hence,  $Q(s)$ is decreasing in every
segment of its continuity.  Make sure that jumps of $Q(s)$
at some point $s_3$, such that
$ k(s_3{-}0) < k(s_3{+}0)$,
conform to its decreasing behavior. Because functions $x(s)$, $y(s)$, 
$\tau(s)$ and $C(s)$ are continuous, and $C(s_3)$ is still negative,
we deduce from~\refeq{Qs}:
\Equa{Qjump}{%
   Q(s_3{+}0)-Q(s_3{-}0)
  = \Frac{1}{4}[k(s_3{+}0) - k(s_3{-}0)]\,C(s_3) < 0,
}
and $Q(s)$ is decreasing in the entire segment $[0,\,s_4]$.

In the case of biarc in $[0,s_1] \cup [s_1,s_2]$,
$k(s)$ is piecewise constant.
Function $Q(s)$ is continuous and zero up to $s_1$,
remains so in $s_1$ (despite of curvature jump,
due to \EQ{C(s_1)}),
and until $s_2$, due to \EQ{\D{k}(s)} in~\refeq{dQs}.
If the entire curve is biarc, the initial circle of curvature\Kl{1} 
is never again reached by
the second subarc until it makes complete $2\pi$-turn,
which contradicts to Def.~\ref{DefSpiral0}.

At $s\igt s_2$, $Q(s)$ either continuously decreases,
or undergoes negative jumps like~\refeq{Qjump}.
The theorem holds in $[0,s_4]$. 
It remains to prove that the point $Z(s_4)$ does not exist.
Under conditions \EQ{C(s_4)} and \LT{Q_4\ieq Q(s_4)} at such point, 
an attempt to determine $Z(s_4)\ieq x_4{+}\iu y_4$ from two equations~\refeqeq{Cs}{Qs}
fails:
excluding $x$ yields
\equa{%
  \begin{array}{l}
    k_1^2y_4^2 - 2k_1y_4(1{+}2Q_4\cos\tau_4 {-}\cos\tau_4) +
    (2Q_4{-}1{+}\cos\tau_4)^2=0,\\
    y_4=\Frac{1}{k_1} \left[
        1-\cos\tau_4(1{-}2Q_4) \pm \sin\tau_4\sqrt{Q_4(1{-}Q_4)}\,
        \right],
  \end{array}
}
i.e. unsolvability with \LT{Q_4} (or immediate contradiction \LT{\sin^2(\tau_4/2)}, 
if \EQ{k_1}).
So, 
spiral never returns to its initial circle of curvature \Kl{1}.
If the curvature decreases, the curve deviates to the
right of \Kl{1}, and $C(s)$ changes sign. 
So do derivative $\D{k}(s)$ and curvature jumps, thus preserving
inequalities~\refeq{Qjump},~\refeq{dQs}, and~\refeq{Qle0}.\qed{}
\end{pf}
%
%
The corollary to this theorem, due to W.~Vogt 
(Satz~1 in~\cite{Vogt}), is the absence of double points
on a spiral. 
Kneser's theorem (see \cite{Guggen}, theorem~3--12), stating that 
{\em``Any circle of curvature of a spiral arc 
contains every smaller circle of curvature of the arc in its interior
and in its turn is contained in the interior of every circle of
curvature of greater radius''} concerns spirals without inflection
and can be generalized as follows:
\begin{cor}
\label{MaterialCor}
Let $\Kl{}(s)=\Kr{x(s),y(s),\tau(s),k(s{\pm}0)}$ be a family of
circles of curvature of a spiral curve. Then the region of material
of any circle includes the region of material of any other circle
with greater curvature:
\equa{
   k(s_2)> k(s_1) \So \Mat{\Kl{}(s_2)}\subset \Mat{\Kl{}(s_1)}.
}
\end{cor}
\Reffig{figure05} illustrates the statement. 
The region $\Mat{\Kl{M}}$ of the initial circle of curvature 
is the whole plane except the interior of \Kl{M}.
As the curvature increases,
the regions of material become smaller, each next being within
the previous one. They remain unbounded up to the inflection point,
whose region of material is the right half-plane.

\Reffigs{figure06}a,b,c illustrate theorem~\ref{MainTheorem} for normalized
spiral. With boundary circles of curvature~\refeq{K1K2c1}
and angles $\omega$ and $\gamma$, defined by~\refeq{OmegaGamma},
inequality~\refeq{Qle0} takes form
\Equa{Qdef1}{%
  \begin{array}{rcll}
     Q(\kn_1, \kn_2,\alpha,\beta)
         &=& \kn_1 \kn_2+ \kn_2\Sa- \kn_1\Sb+\sin^2\gamma &{}={}\\
         &=&( \kn_1+\Sa)\,( \kn_2-\Sb)+\sin^2\omega &{}\le 0,
  \end{array}
}
%
\Lwfig{p}{0.9\textwidth}{figure06}{}
Having fixed~$\alpha$ and~$\beta$, consider the region of permissible values 
for normalized boundary curvatures in the plane $(\kn_1,\kn_2)$. 
This region consists of two subregions,
each bounded by
%
%
one of the two branches of the hyperbola
\EQ{Q(\kn_1, \kn_2)}, traced at the left side of \Reffig{figure06}.
Its centre is located in the point
$C(\kn_1,\kn_2)\ieq(-\sin\alpha,\sin\beta)$;
these two curvatures correspond to those of lense's boundaries.

Biarcs marked as $h_i$ in the right side have boundary curvatures 
$(\kn_1,\kn_2)$ corresponding to the points $H_i$ 
of the hyperbola. By theorem~\ref{MainTheorem},
every biarc represents the {\em unique spiral}, 
matching end conditions of this kind.
Non-biarc curves are presented by some point~$K$
in the plane $(\kn_1,\kn_2)$,
and the arc~$k$ of Cornu spiral in the plane $(x,y)$.
\begin{cor}
\label{k1k2Cor}
End conditions of a normalized spiral arc
obey the following inequalities:
\Equa{k1k2}{%
\begin{array}{rcl}
       \kn_1< {-}\sin\alpha,\quad \kn_2>\sin\beta,
       &\quad  &\mbox{if~~}\kn_1< \kn_2;\\
       \kn_1> {-}\sin\alpha,\quad \kn_2<\sin\beta,
       &\quad  &\mbox{if~~}\kn_1> \kn_2.
\end{array}
}
\end{cor}
\begin{pf}
Inequalities~\refeq{k1k2} merely reflect the position
of two regions \LE{Q} with respect to the asymptotes of the 
hyperbola,
vertical ($\kn_1\ieq{-}\sin\alpha$) and horizontal one ($\kn_2\ieq\sin\beta$).
It remains to show that the line $\kn_2\ieq \kn_1$
separates the two branches of hyperbola, 
thus connecting inequalities~\refeq{k1k2}
to specified conditions, increasing or decreasing curvature.
Substituting $\kn_1\ieq \kn_2\ieq \kn$ into the equation 
\EQ{Q(\kn_1,\kn_2,\alpha,\beta)} yields
\Equa{k1eqk2}{%
  \kn^2+\kn(\Sa{-}\Sb)+\sin^2\gamma =
  (\kn{+}\sin\gamma\cos\omega)^2 +\sin^2\gamma\sin^2\omega = 0.
}
Hence, except special cases \EQ{\sin\gamma} or \EQ{\sin\omega},
statement~\refeq{k1k2} is valid:
two convex regions, bounded by the upper left and
the lower right branches of hyperbola,
are the\linebreak
regions of possible boundary
curvatures for $\kn_1\ilt \kn_2$ and $\kn_1\igt \kn_2$ respectively.

The first exception, \EQ{\sin\gamma}, occurs if
$\alpha\ieq\beta$ and provides the unique common point
\EQ{\kn_1\ieq \kn_2} without intersecting the hyperbola
(point $H_0$ and degenerate biarc $h_0$ in \Reffig{figure06}b).
The statement, assuming $\kn_1\ineq \kn_2$, remains valid.

The second exception, \EQ{\sin\omega}, i.e.  $\alpha\ieq{-}\beta$, 
is illustrated by \Reffig{figure06}d. The hyperbola degenerates into a pair 
of straight lines with the centre $C$ on the line $\kn_1\ieq \kn_2$, 
still separating two regions in question.
However, inequalities~\refeq{k1k2} should be considered as non-strict,
like $\kn_1\ile{-}\sin\alpha$, because the hyperbola is coincident with 
its asymptotes.
By theorem~\ref{MainTheorem}, a~spiral, corresponding to equality, 
may be only biarc. Attempt to construct it gives the only possibility:
the first subarc with $\kn_1\ieq{-}\sin\alpha$, going from~$A$ to~$B$,
and the second subarc being
circumference of a circle of any curvature $\kn_2$ from~$B$ to~$B$
(examples $h_1$ and $h_2$). Similar constructions $h_{3,4}$
arise with $\kn_2\ieq\sin\beta$ and arbitrary $\kn_1$.
Since Def.~\ref{DefSpiral0} excludes this construction,
the points of such degenerated hyperbola do not produce a spiral,
and should be excluded. 
Inequalities~\refeq{k1k2} remain strict.\qed{}
\end{pf}
%
Let us apply inequality~\refeq{Qle0} to another Vogt's statement,
namely, that {\em spiral has no double tangent}
(Satz~7 in~\cite{Vogt}).
Its refined form is illustrated by the double tangent $\Vec{BA}$
in \Reffig{figure05}, and sounds like 
\begin{cor}
Spiral curve may have double tangent only if this tangent
joins points with opposite curvature sign, or is coincident
with the inflection segment of the spiral.
\end{cor}
\begin{pf}
Two curvature elements with common tangent
can be denoted as
\equa{%
  K(x_1,y_1,\tau_1,k_1) \mbox{~~~and~~~}
  K(x_1{+}t\cos\tau_1,\,y_1{+}t\sin\tau_1,\,\tau_1,k_2),\quad t\neq 0.
}
Inequality~\refeq{Qle0} takes form \LE{4Q\ieq k_1 k_2 t^2},
supplying the proof for the general case, \LT{Q},  
\NE{k_{1,2}}. Consider exceptions.
If, say, \EQ{k_1}, the spiral deviates from its tangent as soon 
as $k(s)$ becomes non-zero. By Cor.~\ref{MaterialCor},
the spiral has no more common points with this tangent. 
The biarc case, \EQ{Q}, is trivial: two tangent circles may have
common tangent only in their unique common point.\qed{}
\end{pf}
\begin{cor}
\label{InflectionCor}
The tangent
at the inflection point of a normalized spiral arc
cuts the interior of the chord and is directed downwards
$($i.e. \LT{\sin\tau(s_0)}$)$ if curvature increases,
or upwards $($\GT{\sin\tau(s_0)}$)$ if curvature decreases.
\end{cor}
\begin{pf}
The tangent at the inflection point $Z(s_0)$ is
at the same time the circle of curvature~$\Kl{0}$
(\Reffigs{figure05} and \Figref{figure11}a).
By Cor.~\ref{MaterialCor}, the endpoints $A$ and $B$ of the arc
are disposed bilaterally along \Kl{0}.
That's why \Kl{0} cuts the chord in the interior.
For increasing $k(s)$, point~$A$ is located 
from the right of the line \Kl{0}, and $B$ from the left of it. 
For normalized curve, when $\Vec{AB}$ is brought  horizontal,
this is equivalent to downwards directed tangent~\Kl{0}.
\qed{}
\end{pf}
Two following propositions are our previous results from~\cite{InvInv}.
They can be easily derived from parametric equation of curve, inverse to given one,
by calculating and differentiating its curvature.
\begin{prop}
\label{InversionProp}
Inversion, applied to a spiral curve, preserves the monotonicity 
of the curvature, interchanging its decreasing/increasing character.
\end{prop} 
\begin{prop}
\label{InvCurvatureProp}
If a curvature element~$\Kl{1}$ is inverted with respect to
a circle of inversion~$\Kl{0}$, the curvature of the image~$\Kl{2}$
is given by 
\equa{
      k_2 = 2k_0(1-2Q_{01}) - k_1,\qquad Q_{01} = Q(\Kl{0},\Kl{1}).
}
\end{prop}
The direction, artificially assigned to the circle of inversion, 
does not affect the inverse curve. 
If $\Kl{0}$ is reversed, both $k_0$ and $(1{-}2Q_{01})$ change sign.

\Topic{Vogt's theorem for long spirals}
%
On the spiral $Z(s)\ieq x(s){+}\iu y(s)$, $s\In[0,S]$,
consider subarc $s\In[u,v]$ and define functions
\Equa{Defmu}{%
    h(u,v)=|Z(v)-Z(u)|,\qquad \mu(u,v)=\arg[Z(v)-Z(u)]
}
for the length and direction of the chord.
For any subarc and the entire curve the cumulative boundary angles $\acum(u,v)$
and $\bcum(u,v)$ with respect to varying chord $AB(u,v)$ can be expressed as
\Equa{Defabuv}{%
\begin{array}{l}
  \acum(u,v) = \tau(u)-[\mu(u,v)+2\pi m],\\
  \bcum(u,v) = \tau(v)-[\mu(u,v)+2\pi m],
\end{array}
}
satisfying the natural condition for the {\em winding angle}~$\rho$
of the arc:
\equa{%
     \bcum(u,v) - \acum(u,v) = 
     \tau(v)-\tau(u) = \Int{u}{v}{k(s)}{s} = \rho(u,v).
}
This still allows to assign any value $\alpha{+}2\pi n$ to $\acum$.
To fix it, we note that the angles $\acum$ and $\bcum$
can be unambiguously determined within the range $(-\pi,\pi)$ 
for rather short subarc $[u_0,v_0]$. 
In particular, \ \EQ{\acum(u,u)\ieq\bcum(u,u)}. 
Define cumulative angular functions for any arc $[u,v]$ as
\Equa{ABcum}{%
     \acum,\bcum(u,v) =
     \lim\limits_{{\displaystyle {}^{u_0\to u}_{ v_0\to v} }}
     \acum,\bcum(u_0,v_0)\,,\quad u\le u_0 = v_0 \le v,
}
preserving continuity at $\acum,\bcum=\pm\pi,\,\pm3\pi,\,\ldots\,,$ 
while the limits are being reached.
\begin{lem}
\label{CumLemma}
Cumulative boundary angles $\acum,\,\bcum$, defined by
Eq.~\refeq{ABcum}, do not depend
on the start point $u_0\ieq v_0$ and the way
the limits are reached.
\end{lem}
\begin{pf}
To calculate $\acum(u,v)$ and $\bcum(u,v)$ or,
for symmetry, function
\Equa{OmegaCum}{%
     \ocum(u,v) = \Frac{1}{2}[\acum(u,v)+\bcum(u,v)]
                = \Frac{1}{2}[\tau(u){+}\tau(v)]-[\mu(u,v)+2\pi m],
}
let us restore it from derivatives, which are free from
$2\pi m$-uncertainty:
\equa{%
  \begin{array}{lcl}
      \Pd{\mu(u,v)}{u}&=&
          \Pd{}{u} \arctan\Frac{y(v)-y(u)}{x(v)-x(u)}=\\[8pt]
           &=&\Frac{-\sin\tau(u)\,[\overbrace{x(v){-}x(u)}^{h\cos\mu}]
                 +\cos\tau(u)\,[\overbrace{y(v){-}y(u)}^{h\sin\mu}]}%
                {[x(v){-}x(u)]^2 + [y(v){-}y(u)]^2} 
      = -\Frac{\sin[\overbrace{\tau(u){-}\mu(u,v)}^{\strut\alpha(u,v){+}2\pi m}]}{h(u,v)}.
  \end{array}
}
So,
\equa{%
    \Pd{\mu(u,v)}{u}=-\Frac{\sin\alpha(u,v)}{h(u,v)},\qquad
    \Pd{\mu(u,v)}{v}= \Frac{\sin\beta(u,v)}{h(u,v)},
}
and
\Equa{domega}{
\begin{array}{l}
   \Pd{\ocum(u,v)}{u} =
   \Frac{1}{2}k(u) + \Frac{\sin\alpha(u,v)}{h(u,v)} = G(u,v),\\[10pt]
   \Pd{\ocum(u,v)}{v} =
   \Frac{1}{2}k(v) - \Frac{\sin\beta(u,v)}{h(u,v)} = H(u,v).
\end{array}
}
\begin{figure}[t]
\Bfig{.22\textwidth}{figure07}{}
\hfill  
\Bfig{.66\textwidth}{figure08}{}
\end{figure}
Continuous function $\ocum(u,v)$ can be now calculated as
\Equa{IntGH}{%
     \Int{W_i W}{}{G(u,v)}{u}+H(u,v)\,{\rm d}v
}
along any path $W_i W$ (shown in \Reffig{figure07})
in the closed triangular region $0\ile u\ile v\ile S$
with $W_i$ taken on the boundary $u\ieq v$.
From the expansions
\equa{%
 \begin{array}{l}
   \tau(u{+}s) = \tau(u) + s k(u) +O(s^2),\\[3pt]
   Z(u{+}s) = Z(u) + s\e^{\iu\tau(u)}
                 + \frac{\iu}{2}s^2k(u)\e^{\iu\tau(u)} + O(s^3),\\[2pt]
   h(u,u{+}s) = |Z(u{+}s){-}Z(u)| = s+O(s^2),\\[3pt]
     \mu(u,u{+}s) = \arg[Z(u{+}s)-Z(u)] = 
               \arg[\e^{\iu\tau(u)}(s + \frac{\iu}{2}s^2k(u) + O(s^3))] =\\[2pt]
   \hphantom{\mu(u,u{+}s)} = \arg[\e^{\iu\tau(u)}(1 + \frac{\iu}{2}sk(u) + O(s^2))]
   =
                \tau(u) + \frac{1}{2}sk(u) + O(s^2),\\[3pt]
  \alpha(u,u{+}s) = \tau(u)-\mu(u,u{+}s)=-\frac{1}{2}sk(u) +O(s^2),\\[3pt]
  \beta(u,u{+}s) = \tau(u{+}s)-\mu(u,u{+}s) = \frac{1}{2}sk(u) +O(s^2),
 \end{array} 
}
it follows that derivatives \refeq{domega} can be continuously
defined on the line $u\ieq v$ as zeros:
\equa{
  G(u,u)=\lim\limits_{s\to 0}G(u,u{+}s)
        =\lim\limits_{s\to 0}
        \left[\Frac{1}{2}k(u)-\Frac{\sin[s k(u)/2]}{s}\right]
        = 0,\quad H(u,u)\ieq\ldots\ieq 0
}
This yields \EQ{\int_{W_1W_2}} along the line $u\ieq v$; and, together with
\equa{%
\Pd{G(u,v)}{v} =
\Pd{H(u,v)}{u}\; \left[= -\Frac{\sin(\alpha{+}\beta)}{h}\right],
}
ensures the independence of the integral~\refeq{IntGH}
on a start point $W_i$ and integration path.
In terms of definition~\refeq{ABcum}, the limits
do not depend on the start point
and the way they are reached. Finally, 
$\acum\ieq\ocum{-}\rho/2$,
$\bcum\ieq\ocum{+}\rho/2$.\qed
\end{pf}
\begin{defn}
\rm
The angle $\ocum$, defined above, will be referred to
as {\it Vogt's angle} of a spiral arc.
For short spiral Vogt's angle is the signed half-width of the 
lense~$\Lense{\alpha,-\beta}$.
\end{defn}
\newcommand{\sref}{s_{0}} 
\begin{defn}
\rm
By the {\em reference point\/} of the spiral shall be meant
the point $Z(\sref)$, corresponding to the minimal absolute value
of curvature
(points $A_1$, $B_2$ and $O$ in \Reffig{figure08}).
If the spiral has inflection, $Z(\sref)$ is the inflection point;
otherwise it is either start point ($\sref\ieq0$)
or the end point ($\sref\ieq S$).
\end{defn}
If curvature increases (decreases), function $\tau(s)$
is downwards (upwards) convex and attains minimal (maximal) value
at the reference point.
\begin{lem}
\label{TauPiLemma}
Function $\tcum(s)$ for normalized spiral,
taking cumulative angles $\tcum(0)\ieq\acum$ and $\tcum(S)\ieq\bcum$
as the boundary values, can be distinguished from its other versions,
$\tcum(s)\pm 2\pi m$, by the value at the reference point:
\end{lem}
\Equa{TauCum}{%
  0< |\tcum(\sref)|<\pi,\mbox{~~~~or~~~~}
   \begin{array}{ll}
       -\pi<\tcum(\sref)<0   &\mbox{~~~if~~}k_1< k_2,\\
     \HM{}0<\tcum(\sref)<\pi &\mbox{~~~if~~}k_1> k_2.
   \end{array}
}
\begin{pf}
Based on lemma~\ref{CumLemma}, calculate angles
$\acum(0,S)$ and $\bcum(0,S)$ as follows.
Choose the coordinate system with the origin at the reference point $Z(\sref)$
with X-axis directed along $\tau(\sref)$ (\Reffig{figure08}).
In this system function $\tau(s)$ can be defined as
\equa{%
   \tau(s)=\Int{\sref}{s}{k(\xi)}{\xi},\qquad \tau(\sref)=0.
}
Calculate limits~\refeq{ABcum} starting from $u_0\ieq v_0 \ieq \sref$,
$u$~non\-inc\-rea\-sing,
$v$~non\-de\-crea\-sing.
Consider the case of increasing curvature.
If curvature is nonnegative (spiral $A_1B_1$), then,
by Cor.~\ref{MaterialCor}, the curve is located in the upper half-plane,
except initial subarc of zero curvature, if any.
We keep $u$ constant (\EQ{u\ieq\sref}) and increase $v$. 
Because
\equa{%
 \mu(\sref,\sref)=
 \lim\limits_{\varepsilon\to {+}0}\mu(\sref,\sref{+}\varepsilon)
 =\tau(\sref)=0,
}
the angle $\mu(u,v)$ never attains the values $\pm\pi$,
becoming strictly positive as soon as the point $Z(v)$
deviates from the axis of $X$:
\Equa{MuPi}{%
     0<\mu(u,v)<\pi\qquad
     \mbox{if~~}u\In[0,\sref],\quad v\In[\sref,S],
     \mbox{~~and~~} k(u)< k(v).
}
The case of increasing nonpositive curvature (spiral $A_2B_2$)
is similar: we start with $u_0\ieq v_0\ieq\sref\ieq S$ (point $B_2$),
$u$ decreasing, $v\ieq S$ kept constant. Curve being located in
the lower half-plane, \refeq{MuPi} remains valid.
So it does in the inflection case ($A_3B_3$):  points $Z(u)$, $u\In[0,\sref]$,
are located in the lower half-plane or its boundary,
points $Z(v)$, $v\In[\sref,S]$~--- in the boundary or upper half-plane.
The angle $\mu(u,v)$ 
becomes strictly positive as soon as the point $Z(u)$ or $Z(v)$
deviates from the axis of $X$, and never attains the value $\pi$.
For the three cases $2\pi m$-uncertainty
in~\refeq{Defabuv} disappears, \EQ{m}, and
\equa{%
   \acum = \tau(0)-\mu(0,S),\qquad
   \bcum = \tau(S)-\mu(0,S).
}
To normalize the curve,
chords $A_iB_i$ should be brought to horizontal position,
i.e. rotation through the angle \ $-\mu(0,S)\In(-\pi,0)$ should
be applied. This means replacing $\tau(s)$ by
\equa{%
   \tcum(s) = -\mu(0,S)+\tau(s),
}
whose values at $s=0$,~$S$ and~$\sref$ are equal to
\equa{%
   \tcum(0)=\tau(0)-\mu(0,S)=\acum,\quad
   \tcum(S)=\bcum,\qquad
   \tcum(\sref)= 0-\mu(0,S)\in(-\pi,0).\quad\qed{}
}
\end{pf}
The angle $\mcum(u,v)$ can be defined in the cumulative sense
similarly to $\ocum(u,v)$; limits like~\refeq{ABcum} for~$\mcum$
can be calculated starting from $\mcum(s,s)\ieq\tcum(s)$.
In an equivalent manner $\mcum(u,v)$ can be derived 
from~\refeq{OmegaCum}, where it is bracketed as $[\mu(u,v)+2\pi m]$:
\Equa{MuCum}{%
    \mcum(u,v) = \Frac{\tcum(u){+}\tcum(v)}{2}-\ocum(u,v).
}
For future use we notice, as an immediate corollary of~\refeq{MuPi}, 
the following inequality:
\Equa{DeltaPi}{%
        -\pi < \mcum(\sref,S)-\mcum(0,\sref) < \pi.
}
\begin{thm}
\label{LongVogtTheorem}
With boundary angles, defined in cumulative sense, Vogt's theorem
remains valid for a spiral of any length:
\equa{%
\sign\ocum(u,v) = \sign[k(v){-}k(u)],\quad\mbox{or}\quad
\sign(\acum{+}\bcum) = \sign(k_2{-}k_1). 
}
Except the circular subarc,
wherein Vogt's angle $\ocum(u,v)$ is constant and zero,
it is strictly monotonic function of the arc boundaries:
\equa{%
        |\ocum(u_1,v_1)| < |\ocum(u,v)|
        \mbox{~~~~if~~~}  [u_1,v_1] \subset [u,v],
        \mbox{~~~~and~~~} k(u)\neq k(v).
}
\end{thm}
\begin{pf}
Rewrite derivatives~\refeq{domega} involving normalized curvatures
$\kn_1\ieq \frac{1}{2}h(u,v)k(u)$, 
$\kn_2\ieq \frac{1}{2}h(u,v)k(v)$,
and, assuming increasing curvature,
the first row of~\refeq{k1k2}:
\Equa{dOmgdudv}{%
\Pd{\ocum(u,v)}{u} =
\Frac{\kn_1+\sin\alpha}{h}\le 0,\qquad
\Pd{\ocum(u,v)}{v} =
\Frac{\kn_2-\sin\beta}{h} \ge 0.
}
Equalities are added to account for the possible occurrence 
of a circular subarc within the spiral. If it is not the case,
or as soon as $k(u)\ineq k(v)$, function $\ocum(u,v)$ grows up 
when $u$ decreases and/or $v$ increases.
Because \EQ{\ocum(s,s)}, \ \GT{\ocum(u,v)}  follows.\qed{}
\end{pf}
Recall that Vogt's angle is in fact the intersection angle of 
two circles. Taking into account continuity of $\ocum(u,v)$, and Prop.~\ref{InversionProp},
we conclude the following:  
\begin{cor}
Inversion changes sign of Vogt's angle, preserving its absolute value.
\end{cor}
\Lwfig{t}{0pt}{figure09}%

Now consider the point $P\ieq Z(s)$, moving along the normalized spiral
from $A$ to $B$, and two subarcs of the spiral, $AP$ and $PB$
(\Reffig{figure09}). Denote  
$h_1(s)\ieq h(0,s)\ieq|AP|$ and 
$h_2(s)\ieq h(s,S)\ieq|BP|$,
and apply similar notation to introduce functions
\equa{%
   \begin{array}{llll}
      \alpha_1(s)\ieq \acum(0,s),\quad&
      \beta_1(s) \ieq \bcum(0,s),\quad&
      \omega_1(s)\ieq \ocum(0,s),\quad&
      \mu_1(s)   \ieq \mcum(0,s),\\
      \alpha_2(s)\ieq \acum(s,S),\quad&
      \beta_2(s) \ieq \bcum(s,S),\quad&
      \omega_2(s)\ieq \ocum(s,S),\quad&
      \mu_2(s)   \ieq \mcum(s,S).
   \end{array}
}
From equations
\equa{%
  \begin{array}{lclclc}
    x(s)&{}={}&
    h_1(s)\cos\mu_1(s)-c&{}={}&-h_2(s)\cos\mu_2(s)+c,\\
    y(s)&{}={}& 
    h_1(s)\sin\mu_1(s)  &{}={}&-h_2(s)\sin\mu_2(s)
  \end{array}
}
it follows that
\Equa{BiPolar}{%
h_1(s) = \Frac{2c\sin\mu_2(s)}{\sin\delta(s)},\quad
h_2(s) = \Frac{-2c\sin\mu_1(s)}{\sin\delta(s)},\quad
\mbox{where~~~}\delta(s)=\mu_2(s)-\mu_1(s).
}
Function $\delta(s)$ inherits continuity and cumulative treatment
from $\mu_{1,2}(s)$. It is the turning
of the chord's direction at~$P$,
or signed external angle of the triangle $APB$ at~$P$.
If $P$ is in the upper half-plane,  $\sin\delta(s)$
is negative; it is positive in the lower half-plane;
$\delta(s)\ieq 2\pi n$ if $P$ is within the chord,
and $\delta(s)\ieq\pi(2n{-}1)$ if $P$ belongs to chord's complement.
The locus of points, where $\delta$ is constant,
is the arc  $\Aarc{-\delta}$.
\begin{lem}
\label{DeltaLemma}
Function $\delta(s)$, defined on a spiral with 
increasing $(decreasing)$ curvature, 
is strictly increasing $(decreasing)$ from
$\delta(0)\ieq {-}\acum$ to $\delta(S)\ieq\bcum$,
taking value $-\pi\ilt\delta(\sref)\ilt\pi$ at the reference point;
its derivative is continuous and does not vanish in $[0,S]$.
\end{lem} 
\begin{pf}
Using \refeq{MuCum}, rewrite $\delta(s)=\mu_2(s){-}\mu_1(s)$ as
\Equa{DeltaCum}{
  \delta(s) = 
        \left[\Frac{\tcum(s){+}\bcum}{2}-\ocum(s,S)\right] -
        \left[\Frac{\acum{+}\tcum(s)}{2}-\ocum(0,s)\right] =
         \Frac{1}{2}\rho(0,S)+\omega_1(s)- \omega_2(s).
}
Apply~\refeq{domega} to calculate derivative:
\begin{eqnarray}
    \D{\delta}(s)
       &=& \D{\ocum}_1(s)-\D{\ocum}_2(s) =
           \Pd{\ocum(u,v)}{v}\Big|^{v=s}_{u=0} -
           \Pd{\ocum(u,v)}{u}\Big|^{v=S}_{u=s} =\nonumber\\[8pt]
       &=&
           \left[\Frac{1}{2}k(s)-\Frac{\sin\beta_1(s)}{h_1(s)}\right]-
           \left[\Frac{1}{2}k(s)+\Frac{\sin\alpha_2(s)}{h_2(s)}\right]
           = 
           -\Frac{\sin\beta_1(s)}{h_1(s)}-
            \Frac{\sin\alpha_2(s)}{h_2(s)}.\label{DeltaInc}
\end{eqnarray}
Because $k(s)$ disappears in~\refeq{DeltaInc}, 
$\D{\delta}(s)$ is continuous even if the curvature jump occurs
(smooth plot $\omega_1{-}\omega_2$ in \Reffig{figure09},
compared to $\omega_1$ and $\omega_2$, illustrates it).

As the point $P(s)$ moves along the spiral, the arc $AP$
is lengthening, and the arc $PB$ shortening.
From theorem~\ref{LongVogtTheorem}
it follows that $\omega_1(s)$ is increasing,
and $\omega_2(s)$  decreasing (the case of increasing curvature
is being considered). The difference $\omega_1(s){- }\omega_2(s)$
in~\refeq{DeltaCum} is therefore an increasing function.
So is $\delta(s)$: \ \GE{\D{\delta}(s)}. 
The only possibility for equality is that $\D{\ocum}_{1,2}(s)$ 
are simultaneously zeros, i.e. subarcs $AP$ and $PB$ are of constant 
curvature. This means that the spiral $APB$ is biarc, 
and $P$ is the point of tangency of its two circular arcs,
as depicted in \Reffig{figure09}. 
If it is the case, then the last two fractions in~\refeq{DeltaInc}
are equal to $k_1$ and $-k_2$ respectively, and
\GT{\D{\delta}(s_P)\ieq\frac{1}{2}(k_2{-}k_1)}.
The derivative is thus strictly positive everywhere in $(0,S)$,
$\delta(s)$ is strictly increasing.

To calculate derivatives at the endpoints we have to resolve uncertanties $0/0$ 
in~\refeq{DeltaInc}. For $\D{\delta}(0)$ approximate the spiral near the startpoint
by its circle of curvature:  
%
%
\equa{%
\begin{array}{lcl}
  Z(s)& =& -c + \Frac{\iu}{k_1} \e^{\iu\alpha}(1- \e^{\iu k_1 s}) =  
         -c +\e^{\iu\alpha}s + \Frac{\iu}{2}\e^{\iu\alpha} k_1s^2+O(s^3),\\[4pt]
  \delta(s)&=&\arg[c-Z(s)] - \arg[Z(s)+c]=\\[6pt]
           &=&\left[-\Frac{\sin\alpha}{2c}s+O(s^2)\right] -
            \left[\alpha+\Frac{1}{2}k_1 s+O(s^2)\right]=
            -\alpha-\Frac{k_1c{+}\sin\alpha}{2c} s+O(s^2).         
\end{array}
}
The coefficient at $s$ is the derivative $\D{\delta}(0)$. 
Similarly $\D{\delta}(S)$ at the endpoint~$B$ can be found. 
Both are positive due to~\refeq{k1k2}:
\equa{
     \D{\delta}(0)= -\Frac{k_1c{+}\sin\alpha}{2c}>0,\qquad 
     \D{\delta}(S)= \Frac{k_2c{+}\sin\beta}{2c}>0. 
} 
The value of $\delta(\sref)$ at the reference point 
is already estimated in~\refeq{DeltaPi}.
Boundary values can be calculated from~\refeq{DeltaCum}:
\equa{%
  \delta(0) = \Frac{\bcum{-}\acum}{2} + 0 -  \Frac{\acum{+}\bcum}{2}
            = {-}\acum,\qquad
  \delta(S) = \Frac{\bcum{-}\acum}{2} + \Frac{\acum{+}\bcum}{2} - 0
           = \bcum. \;\qed
}
\end{pf}
\Lwfig{t}{0pt}{figure10}%
The plots $\omega_1{+}\omega_2$ in \Reffigs{figure09} and~\Figref{figure10}
illustrate the following property of Vogt's angles: however large be the
range $[0,\ocum]$ of the monotonic functions $\omega_1(s)$ and $\omega_2(s)$,
their sum is enclosed in the interval of the width $\pi$:
\begin{lem}
Let \ $\Omega(s) = |\omega_1(s)|+|\omega_2(s)|-|\ocum|$. Then
\Equa{o1o2}{%
     -\pi < \Omega(s) < 0 \mbox{~~~for~~~}  s\in(0,S).
}
\end{lem}
\begin{pf}
For the case of increasing curvature all $\ocum$'s are nonnegative,
and
\equa{%
  \begin{array}{rcl}
     \Omega(s) &=&
     \omega_1(s)+\omega_2(s)-\ocum \eqref{OmegaCum}
     \left[\Frac{\acum{+}\tcum(s)}{2}-\mu_1(s)\right] +
     \left[\Frac{\tcum(s){+}\bcum}{2}-\mu_2(s)\right] -
          \Frac{\acum{+}\bcum}{2} =\\
    &=& \tcum(s)-\mu_1(s)-\mu_2(s).
  \end{array}
}
Continue calculation of the derivative~\refeq{DeltaInc},
which is strictly positive:
\equa{%
  \begin{array}{rcl}
 \D{\delta}(s) &=&
   -\Frac{h_2\sin\beta_1+h_1\sin\alpha_2}{h_1h_2}
   \,\eqref{BiPolar}\,
   2c\Frac{\sin\mu_1\sin\beta_1-\sin\mu_2\sin\alpha_2}{h_1h_2\sin\delta}=\\[8pt]
  &=&
   2c\Frac{\sin\mu_1\sin(\tcum{-}\mu_1)-\sin\mu_2\sin(\tcum{-}\mu_2)}%
      {h_1h_2\sin\delta}
   = \Frac{-2c}{h_1h_2}\sin\Omega(s) > 0.
  \end{array}
}
The function $\Omega(s)$ takes zero values at the endpoints,
and $\sin\Omega(s)$ is strictly negative in $(0,S)$.
Invoking continuity of $\Omega(s)$,
we conclude that it never reaches values 0 or $-\pi$ 
within $(0,S)$. \qed
\end{pf}
Now we recall Def.~\ref{DefShort1} to introduce some
quantitative measure to the notion of a long spiral.
The inflection point, if present,
subdivides spiral into two branches,
left and right. Points $s^{-}_i$, $i \ieq 1,\ldots, M_1$ on the left
branch, and $s^{+}_i$, $i \ieq 1,\ldots, M_2$ on the right branch
are those, where Def.~\ref{DefShort1} of short arc is violated,
i.e. $0\ilt s^{\pm}_i\ilt S$, $\cos\tau(s^{\pm}_i)\ieq{-}1$.
Because \EQ{\sin\tau(s^{\pm}_i)}, this cannot happen at the inflection
(Cor.~\ref{InflectionCor});
so, $s^{\pm}_i$ are distinct points,
not continuous segments, as the inflection could be.
If there are no such points or one branch is absent,
the corresponding  counter~$M_{1,2}$ is zero.

Two other counters, $N_{1,2}$, are introduced
in the context of Def.~\ref{DefShort2}.
They count internal points where the spiral meets
the left ($N_1$) and the right ($N_2$) complements of the chord.

\begin{thm}
\label{MN_Theorem}
Counters $M_{1,2}$ and $N_{1,2}$ are pairwise equal.
Cumulative boundary angles $\acum$ and $\bcum$
for a normalized spiral arc of
increasing/decreasing curvature are
\equa{
   \begin{array}{lllll}
     k_1< k_2:\quad&   \acum\ieq\alpha{+}2\pi N_1,\quad
                     &   \bcum\ieq\beta{+} 2\pi N_2 \quad
                     &  \mbox{\rm~with~~~}
                     &   \alpha,\beta\in(-\pi,\pi],\\
     k_1\igt k_2:\quad&   \acum\ieq\alpha{-}2\pi N_1,\quad
                     &   \bcum\ieq\beta{-} 2\pi N_2 \quad
                     &  \mbox{\rm~with~~~}
                     &   \alpha,\beta\in[-\pi,\pi),
   \end{array}
}
or, rewritten for the case of increasing curvature in more detailed form,
\settowidth{\tmplength}{$\displaystyle{\acum\ieq\alpha\;(N_1\ieq 0),W}$}%
\renewcommand{\tmp}[1]{\makebox[\tmplength][l]{\mbox{$\displaystyle#1$}}}%
\begin{eqnarray}
    0\le k_1< k_2: \quad &
    -\pi<\alpha< 0,\;  -\pi<\beta\le\pi,\quad &
    \tmp{\acum\ieq\alpha\;(N_1\ieq0),}
    \tmp{\bcum\ieq\beta{+}2\pi N_2;}
  \label{CaseInc1}\\
     k_1< 0<\,k_2: \quad &
     -\pi<\alpha\le\pi,\; -\pi<\beta \le\pi,\quad &
     \tmp{\acum\ieq\alpha{+}2\pi N_1,}
     \tmp{\bcum\ieq\beta{+}2\pi N_2;}
  \label{CaseInc2}\\
     k_1< k_2\le 0: \quad &
    -\pi<\alpha\le\pi,\; -\pi<\beta < 0,\quad &
     \tmp{\acum\ieq\alpha{+}2\pi N_1,}
     \tmp{\bcum\ieq\beta\;(N_2\ieq 0).}
  \label{CaseInc3}
\end{eqnarray}
\end{thm}
\Lwfig{t}{0pt}{figure11}%
\begin{pf}
Two cases, \refeq{CaseInc2} and~\refeq{CaseInc1}, are illustrated
by \Reffig{figure11}. Right plots show 
functions $\tcum(s)$ and $\delta(s)$.
Consider the inflection case~\refeq{CaseInc2}, \Reffig{figure11}a.
Function $\tcum(s)$ is decreasing from $\acum$ to its
minimal value $\tcum(\sref)\In(-\pi,0)$
at the reference (inflection) point, and then increases to $\bcum$.
In doing so, it meets $M_1{+}M_2$ times the levels $\pi(2m{-}1)$
(points $T_i$).
The following sequence of its values can be derived
according to definition of counters $M_{1,2}$:
\equa{%
 \acum,\,
 \underbrace{\pi(2M_1{-}1),\, \pi(2M_1{-}3),\,\ldots,\,\pi}_{M_1},\,
 \tcum(\sref),\,
 \underbrace{\pi,\,3\pi,\,\ldots,\,\pi(2M_2{-}1)}_{M_2},\,
 \bcum.
}
By lemma~\ref{DeltaLemma},
function $\delta(s)$ is monotonic increasing from $-\acum$ to $\bcum$.
Spiral cuts the complement of the chord at points $C_i$,
when $\delta(s)$ meets levels $\pi(2n{-}1)$. 
The sequence of its values,
similar to that of $\tcum(s)$, looks like
\equa{%
 -\acum,\,
 \underbrace{-\pi(2M_1{-}1),\, -\pi(2M_1{-}3),\,\ldots,\,-\pi}_{N_1},\,
 \delta(\sref),\,
 \underbrace{\pi,\,3\pi,\,\ldots,\,\pi(2M_2{-}1)}_{N_2},\,
 \bcum.
}
The number of underbraced points is $N_1{+}N_2$ by definition of
counters $N_{1,2}$.
Separating them into two groups is justified as follows.
First, the tangent at the inflection point separates two branches of the
spiral and, by Cor.~\ref{InflectionCor},
two complements of the chord; hence, all the $N_1$ intersections
with the left complement of the chord belong to the left branch
and are followed by the set of the right-sided ones.
Second, by lemma~\ref{DeltaLemma}, the value $\delta(s)\ieq {-}\pi$
terminates the first group of points,
whose number is $N_1$ by definition and $M_1$ by calculation.
Similarly, $M_2\ieq N_2$. 

From the above sequences it follows as well:
if the directions of tangents are $\alpha$ and~$\beta$,
and $|\alpha|,\,|\beta|<\pi$, 
the values $\alpha{+}2\pi N_1$ and $\beta{+}2\pi N_2$ are to be
assigned to cumulative angles $\acum,\,\bcum$.
If $\alpha$ or $\beta$ is equal to $\pm\pi$, the correspondence
is kept by the alternative resolved in favour of $+\pi$ 
(and $-\pi$ in the case of decreasing curvature).

In the case \refeq{CaseInc1} of increasing nonnegative curvature
(\Reffig{figure11}b),
the curve has no left branch, and \EQ{M_1} by definition.
Point~$A$ is the reference point, so, by lemma~\ref{TauPiLemma}, 
$-\pi\ilt\acum\ilt 0$, \ $\acum\ieq\alpha$;
the tangent at $A$ is directed
downwards. The region $\Mat{\Kl{1}}$ is either
half-plane to the left of \Kl{1} (if \EQ{k_1}),
or the interior of the circle \Kl{1} located in this half-plane.
It covers the entire curve (Cor.~\ref{MaterialCor}),
and cannot include any point of the left complement of the chord.
Therefore $N_1\ieq 0\ieq M_1$.

The rest of the proof is similar to that of the inflection case.
Equality $N_2\ieq M_2$ results from monotonic increasing behavior
of functions $\tcum(s)$ (from $\acum$ to $\bcum$) 
and $\delta(s)$ (from $-\acum$ to $\bcum$),
and counting points $T_i$ and $C_i$.

The proofs for the case~\refeq{CaseInc3} of increasing nonpositive curvature
can be obtained by applying symmetry about Y-axis.
Symmetry about X-axis provides the proof for the three similar cases 
of decreasing curvature.\qed{}
\end{pf}
\begin{cor}
\label{EquivDefCor}
Definitions {\,\rm\ref{DefShort1}\,} and {\,\rm\ref{DefShort2}\,} 
of a short spiral are equivalent.
\end{cor}
Polygonal line $ACDF$ in \Reffig{figure12} bounds the open region of
possible values of $\acum,\bcum$ for spirals with
increasing curvature.
Boundary $CD$ results from Vogt's theorem,
$AC$ and $DF$~--- from~\refeq{CaseInc1}--\refeq{CaseInc3}.
Triangle $GDC$, including half-open segments $(CG]$ and $[GD)$,
is the boundary for short spirals.
Biarc curves can be constructed with $\acum,\bcum$ in the interior 
of tra\-pe\-zium $BCDEB$ (see discussion in the next section).
Similar regions for decreasing curvature are symmetric
about the line $CD$ (\EQ{\acum{+}\bcum}). Describing these regions
in the coordinate system $(\rho\ieq\acum{-}\bcum,\,2\ocum\ieq\acum{+}\bcum)$
looks as follows:
\begin{cor}
The winding angle~$\rho$ of a spiral is limited by
\equa{%
      \hphantom{0 < 2|\ocum| <{}} |\rho| < 2|\ocum|+ 2\pi.
}                
If spiral has no inflection {\rm(}cases \refeq{CaseInc1} and \refeq{CaseInc3}{\rm)}, than 
\equa{%
              0 < 2|\ocum| < |\rho| < 2|\ocum|+ 2\pi.
}
\end{cor}
If a spiral undergoes inversion, $|\ocum|$ remains constant, and $\rho$ can vary
in these limits. The latter inequality is similar to the fact 
that a small non-closed circular arc (\EQ{|\ocum|})
can be transformed by inversion to almost $2\pi$-circle, and vise versa
($0\ile |\rho| \ilt 2|\ocum| {+} 2\pi\ieq 2\pi$).

\begin{figure}[t]
\Bfig{.3\textwidth}{figure12}{}
\hfill
\Bfig{.66\textwidth}{figure13}{}
\end{figure}

\Topic{Biarc curves}
%
Biarc curves, considered hitherto
as a flexible tool for curves interpolation,
play an important role in the theory of spiral curves.
However much had been written about biarcs
(see~\cite{Biarc92} and references herein,
\cite{LongBiarcs}~for long biarcs),
the presented description seems to have several advantages.
Normalized position can be considered as canonical
for these curves, and allows to separate the parameters of shape
from positional ones. The proposed parametrization
yields a set of simple and symmetric reference formulae.
No different treatment for ``C-shaped'' and ``S-shaped'' biarcs
is needed.
The specific cases of $\alpha\ieq{\pm}\pi$ or $\beta\ieq{\pm}\pi$,
usually omitted,
are taken into consideration.

The condition of tangency of two arcs, forming biarc,
is the equation
\equa{
   Q(\kn_1, \kn_2, \alpha, \beta)=
   ( \kn_1+\Sa)\,( \kn_2-\Sb)+\sin^2\omega = 0
}
(recall hyperbolas in \Reffig{figure06}). This condition
allows two arcs to be a pair of equally directed straight lines
(``biarc'' $h_0$ in \Reffig{figure06}b).
The hyperbola can be parametrized as follows:
\Equa{BFamily}{%
 \left\{
 \begin{array}{l}
   \kn_1(b)=    -\sin\alpha -  b^{-1}\sin\omega,\\
   \kn_2(b)=\HM{}\sin\beta +  b\sin\omega,
 \end{array}
 \right.
}
Note that $\omega$ is the half-width of the lense,
equal to Vogt's angle $\ocum$ only if biarcs is short;
otherwise $\ocum = \omega \pm \pi$.
Parametrization~\refeq{BFamily}
supplies the parameter $b$ for the one-parametric family
of biarcs with fixed chord $[-1,1]$
and fixed tangent directions $\alpha$ and~$\beta$.
As established in the proof of Cor.~\ref{k1k2Cor}
(\Reffig{figure06}d),
biarcs with \EQ{\sin\omega} do not exist.
Every value of~$b$ produces the unique
point on the hyperbola, unique pair of circles \Kl{1} and \Kl{2},
tangent at the point~$T$, unique path $ATB$,
and unique biarc, denoted below as $\Biarcab{b}{\alpha}{\beta}$
or simply $\Biarc(b)$.

Solving the system of two equations \EQ{C(x,y;\Kl{1,2})}~\refeq{Cxy}
yields coordinates of the point $T\ieq(x_0,y_0)$ of contact of two arcs:
\Equa{PointT}{%
x_0 = \Frac{b^2-1}{\Delta},\quad
y_0 = \Frac{2b\sin\gamma}{\Delta},\qquad
  \Delta = b^2+2b\cos\gamma+1\,.
}
The direction $\tau_0$ of the common tangent at this point is given by
\Equa{Tau0}{%
 \begin{array}{l}  
    \sin\tau_0 \ieq  {-}(b^2\sin\alpha + 2b\sin\omega + \sin\beta)/\Delta,\\
    \cos\tau_0 \ieq\HM{}(b^2\cos\alpha + 2b\cos\omega + \cos\beta)/\Delta,
 \end{array}
 \quad
    \tan\Frac{\tau_0}{2}={-}\Frac%
    {b\sin(\alpha/2)+\sin(\beta/2)}
    {b\cos(\alpha/2)+\cos(\beta/2)}\,.
}
The circular arc from $B$ to $A$, complementary to the arc~$\Aarc{\gamma}$,
will be referred to as the {\em complement of the bisector}
of the lense.
Both bisector and its complement form the circle~$\Gamma$,
shown dashed in \Reffig{figure13}:
\equa{
  \Gamma = K(-1,0,\gamma,{-}\sin\gamma),\qquad
  C(x,y;\Gamma)\eqref{Cxy}{-}\sin\gamma(x^2+y^2-1)-2y\cos\gamma.
}
Except property~\refeq{PropMonoLength}, 
the following properties of biarcs
are either known or easy to prove by means of elementary geometry:
\def\labelenumi{(\theenumi)}        
\def\theenumi{\roman{enumi}}
\begin{enumerate}
\item
\label{PropCircle}
{\em
The locus of contact points $T(b)\ieq(x_0(b),\,y_0(b))$~\refeq{PointT}
is the circle~$\Gamma$. Points $T(b)$ with \GT{b} are located on the 
bisector, those with \LT{b}, on its complement.
}
\smallskip
\item
\label{PropOmega}
{\em
All biarcs meet the the circle~$\Gamma$
at the constant angle~$\omega$.
}
\item
{\em
Definitions
  $\Biarcab{\infty}{\alpha}{\beta}\ieq\Aarc{\alpha}$ and
  $\Biarcab{0}{\alpha}{\beta}\ieq\Aarc{-\beta}$,
are illustrated by \Reffig{figure13} and justified as follows\,}:
\Equa{BiarcLim}{%
 \begin{array}{lllll}
   b\to \infty: &
      \kn_1\to{-}\Sa, \;& \kn_2\to\infty, &\; T\to B,&
           \Biarcab{b}{\alpha}{\beta}\to \Aarc{\alpha};\\
   b\to 0:&
      \kn_1\to\infty,   & \kn_2\to\Sb,\; &\; T\to A,\:&
           \Biarcab{b}{\alpha}{\beta}\to \Aarc{-\beta}.
 \end{array}
}
\item
\label{PropCum}
{\em
The possible values for the pair of counters $(N_1,N_2)$ 
are $(0,0)$, $(0,1)$ or $(1,0)$.
For given tangents $\alpha$, $\beta$,
the cumulative angles
$\acum$, $\bcum$ can take values
{\rm(}arrows mark the cases of increasing $\,({}^\Uparrow)$
or decreasing $\,({}^\Downarrow)$ curvature{\rm):}}
\Equa{BiarcCum}{%
  \begin{array}{ll}
    \mbox{if~~~}\alpha{+}\beta>0,\quad&
    (\acum,\bcum) = (\alpha,\:\beta)^\Uparrow,\quad
                    (\alpha{-}2\pi,\:\beta)^\Downarrow,\quad
                    (\alpha,\:\beta{-}2\pi)^\Downarrow;\\
    \mbox{if~~~}\alpha{+}\beta<0,\quad&
    (\acum,\bcum) = (\alpha,\:\beta)^\Downarrow,\quad
                    (\alpha{+}2\pi,\:\beta)^\Uparrow,\quad
                    (\alpha,\:\beta{+}2\pi)^\Uparrow .
  \end{array}
}
{\em
The first group, biarcs with $(\acum,\bcum)\ieq(\alpha,\beta)$,
corresponds to \GT{b};
they are short and enclosed within the lense.
Biarcs with \LT{b} are located outside the lense.
They are long, unless
\ $\alpha\ieq{\pm\pi}$ or  $\beta\ieq{\pm\pi}$ 
{\rm(see~\ref{PropPi})}.
}
%
\item
\label{PropBinf}
{\em
Discontinuous biarcs {\rm\,(such as shown in
\Reffigs{figure06}a,b by dotted lines)\,}
correspond to nonpositive parameter value~$\binf$, defined by
\equa{
  \binf\ieq
  \left\{
   \begin{array}{lll}
      -\Frac{\sin\omega}{\sin\alpha}\quad&
      \mbox{if~~~} |\alpha| \ige |\beta|\quad & 
      [\,\kn_1(\binf)\ieq 0\,], \\[8pt]
      -\Frac{\sin\beta}{\sin\omega}\quad&
      \mbox{if~~~} |\alpha| \ile |\beta|\quad &  
      [\,\kn_2(\binf)\ieq0\,],
   \end{array} \right. .
}
Three biarcs with $b\In \{\infty;\binf;0\}$ subdivide the XY-plane
into three regions, one of them being the lense. Every region
encloses one of the three subfamilies~\refeq{BiarcCum}. 
}
%
\item
\label{PropPi}
{\em 
If $\alpha\ieq{\pm\pi}$ or  $\beta\ieq{\pm\pi}$ 
{\rm\,(\Reffig{figure06}c)},
one of these three regions, as well as 
one of three subfamilies~\refeq{BiarcCum}, disappear.
The degenerate biarc is at the same time lense's boundary
($\binf\ieq 0$ or $\binf\ieq\infty$).
All such biarcs are short; and
\ $\Biarcab{b}{\pi}{\beta}\ieq\Biarcab{-b}{-\pi}{\beta}$,
\ $\Biarcab{b}{\alpha}{\pi}\ieq\Biarcab{-b}{\alpha}{-\pi}$.
}
%
\item
\label{PropXY}
{\em Taking into account biarcs
$\Biarc(\infty)$, $\Biarc(\binf)$ and $\Biarc(0)$,
there is a unique biarc $\Biarc(b(x,y)\,)$,
passing through every point $(x,y)$ in the plane,
excluding poles $A$ and $B$. Namely,
\Equa{Bxy}{%
   b(x,y) =
   \left\{
   \begin{array}{ll}
     \Frac{\sin\omega[(x{+}1)^2+y^2]}%
          {(1{-}x^2{-}y^2)\sin\alpha-2y\cos\alpha},\quad&
          \mbox{if~~} C(x,y;\Gamma)\sin\omega \le 0,\\[12pt]
     \Frac{(1{-}x^2{-}y^2)\sin\beta+2y\cos\beta}%
          {\sin\omega[(x{-}1)^2+y^2]},\quad&
          \mbox{if~~}  C(x,y;\Gamma)\sin\omega\ge 0.
   \end{array}
   \right.
}
}
%
\item
\label{PropMonoLength}
{\em
The length $L(b)$ of short biarc 
\equa{%
  L(b)=\Frac{\tau_0(b){-}\alpha}{k_1(b)}+
       \Frac{\beta{-}\tau_0(b)}{k_2(b)}\qquad
  \left[ L(0)=\Frac{2c\alpha}{\sin\alpha},\quad
         L(\infty)=\Frac{2c\beta}{\sin\beta}\right],
}
is strictly monotonic function of $b$, 
or constant, if $\alpha\ieq\beta$
{\rm(}the uncertanties $0/0$, wherever occur, are simply reducible{\rm)}.
}
\end{enumerate}
To prove~\refeq{PropCircle} note that~\refeq{PointT} is the parametric
equation of the circle~$\Gamma$. Ordinates of the points, belonging to 
the bisector, should have the same sign as the angle $\gamma$,
and this is achieved at \GT{b}:
\equa{%
    \gamma \neq 0,\quad b > 0  \So  \sign y_0 = \sign \gamma
         \qquad\mbox{(\Reffig{figure06}a)}.
}
If \EQ{\gamma}, the bisector is coincident with the chord;
this corresponds to \GT{b} as well:
\equa{%
    \gamma = 0, \quad b > 0  \So y_0 = 0,\quad 
                 |x_0|=\Big|\Frac{b^2-1}{b^2+2b+1}\Big|=
                 \Big|\Frac{b{-}1}{b{+}1}\Big|  < 1
         \qquad\mbox{(\Reffig{figure06}b)}.
}
In both cases points $T(b)$ for biarcs with \LT{b} fill
the complement of the bisector.

Provided that the locus $T(b)$ is circle~$\Gamma$, 
\refeq{PropOmega} is evident.

Property \refeq{PropCum} can be derived from the fact that 
one of two subarcs of a biarc is located
inside the circle~$\Gamma$, and the other outside it.
Only outside one can meet the chord's complement;
and this may happen only once.

Points of contact $T(b)$, \GT{b}, are inside the lense.
So are subarcs $AT$, $TB$, and the entire biarc curve.
Being inside lense,
{\em biarcs with \GT{b} are short}, their cumulative boundary angles
are $(\acum,\bcum)\ieq(\alpha,\beta)$.
Points $T(b)$, \LT{b}, filling the complement of the bisector,
are outside the lense together with associated biarcs.

The discontinuous biarc $\Biarc(\binf)$ may arise only if one of the
two curvatures is zero. Let \EQ{\kn_1(\binf)} (\Reffig{figure06}a), i.e.
$\binf\ieq{-}\sin\omega/\sin\alpha$. The parametric equation
of the subarc $AT$ is
\equa{%
x(s)=-1+s\cos\alpha,\qquad y(s)=s\sin\alpha.
}
This ray reaches the point of contact $T$ \ \refeq{PointT}
when $y(s_T)\ieq s_T\sin\alpha\ieq y_0$, i.e.
\equa{%
s_T 
    = \Frac{2\binf\sin\gamma}{({\binf}^2{+}2\binf\cos\gamma{+}1)\sin\alpha}
    = \Frac{-2\sin\gamma\sin\omega}%
          {\sin^2\alpha{-}2\cos\gamma\sin\alpha\sin\omega{+}\sin^2\omega}
    = {-}2\Frac{\sin\omega}{\sin\gamma}
}
($\omega{+}\gamma$ was substituted for $\alpha$).
If \GT{s_T}, then tangency occurs before the ray goes off
into infinity, and the biarc is normally continued
by the second arc $TB$ with \NE{\kn_2}. Discontinuity
occurs under the condition $-\infty\ile s_T \ilt 0$,
equivalent to $\cos\alpha\ile\cos\beta$, and
rewritten in~\refeq{PropBinf} as $|\alpha|\ige|\beta|$.

To derive~\refeq{Bxy}, one should first decide, to which subarc
of the sought for biarc $\Biarcab{b}{\alpha}{\beta}$
the point $(x,y)$ belongs.
The circle~$\Gamma$ separates two subarcs, and the decision
depends on the sign of $C(x,y;\Gamma)$. The arc $AT$ goes to the left
of~$\Gamma$ ($AT\In\Mat{\Gamma}$, \ \LE{C(x,y;\Gamma)})
if the vector, defined by $\alpha$,
points to the left of the vector, defined by $\gamma$:
$\gamma<\alpha\ieq\omega{+}\gamma$, i.e. \GT{\omega}, \ \GT{\sin\omega}.
And $AT$ goes to the right of~$\Gamma$ if \LT{\sin\omega}. So,
\equa{%
  \begin{array}{rcl}
    (x,y)\In \Arc{AT}\In\Kl{1}
    &\Longleftrightarrow&
      [C(x,y;\Gamma)\ile0\:\wedge\:\sin\omega\igt0] \;\vee\;
      [C(x,y;\Gamma)\ige0\:\wedge\:\sin\omega\ilt0]\\
    &\Longleftrightarrow&
      C(x,y;\Gamma){\cdot}\sin\omega\ile0.
  \end{array}
}
Under this condition define $b$ from the implicit 
equation~\refeq{Cxy} of circle~\Kl{1}:
\equa{%
     (\underbrace{-\sin\alpha{-} b^{-1}\sin\omega}_{\kn_1})
      \left[(x{+}1)^2{+}y^2\right]+
       2(x{+}1)\sin\alpha - 2y\cos\alpha = 0.
}
%
Below we prove property \refeq{PropMonoLength}.
\begin{pf}
For the case of symmetric lense, i.e. $\alpha\ieq\beta\ieq\omega$,
\EQ{\gamma}, $\tau_0(b)\ieq {-}\omega\ieq const$, and
\equa{%
L(b)=\Frac{-2\omega}{-\sin\omega- b^{-1}\sin\omega}+
       \Frac{2\omega}{\sin\omega+b\sin\omega}
    =\Frac{2\omega}{\sin\omega}\,
}
($c\ieq 1$ assumed).
For general case, \NE{\gamma}, we replace the parameter $b$ by
$\theta\ieq \tau_0(b)$. As $b$ varies from 0 to $\infty$,
$\theta$ varies (monotonously) from $\beta$ to $-\alpha$.
Solving~\refeq{Tau0} for $b$ yields 
\newcommand{\ffrac}[1]{\frac{#1}{2}}
\newcommand{\Sam}[1]{\sin^{#1}\ffrac{\alpha{-}\theta}}
\newcommand{\Sap}[1]{\sin^{#1}\ffrac{\alpha{+}\theta}}
\newcommand{\Sbm}[1]{\sin^{#1}\ffrac{\beta{-}\theta}}
\newcommand{\Sbp}[1]{\sin^{#1}\ffrac{\beta{+}\theta}}
\equa{%
   b=-\Frac{\tan\ffrac{\theta}\cos\ffrac{\beta} +\sin\ffrac{\beta}}%
           {\tan\ffrac{\theta}\cos\ffrac{\alpha}+\sin\ffrac{\alpha}}
    =-\Frac{\Sbp{}}{\Sap{}}\:,
}
and
\equa{%
   \begin{array}{lll}
     \kn_1(\theta)
                   =\sin\gamma\, \Frac{\Sam{}}{\Sbp{}},\;&
     \kn_2(\theta) =\sin\gamma\, \Frac{\Sbm{}}{\Sap{}},\;&
     \kn_2{-}\kn_1=-\Frac{\sin^2\gamma\sin\omega}{\Sap{}\Sbp{}},\\[9pt]
     \D{\kn_1}(\theta)=\Frac{-\sin\gamma\sin\omega}{2\Sbp{2}},&
     \D{\kn_2}(\theta)=\Frac{-\sin\gamma\sin\omega}{2\Sap{2}}.&
   \end{array}
}
The length of biarc as a function of $\theta$, 
and its derivative appear as
\equa{%
  \begin{array}{lcl}
     L(\theta)&=&\Frac{\theta{-}\alpha}{\kn_1(\theta)}+
          \Frac{\beta{-}\theta}{\kn_2(\theta)}\,;     \\[9pt]
     \D{L}(\theta)&=&\Frac{ \kn_1-(\theta{-}\alpha)\D{\kn_1}}{\kn_1^2}-
                  \Frac{\kn_2+(\beta{-}\theta)\D{\kn_2}}{\kn_2^2} =
     \Frac{\kn_1\kn_2(\kn_2{-}\kn_1)-
           \D{\kn_1}\kn_2^2(\theta{-}\alpha)-
           \D{\kn_2}\kn_1^2(\beta{-}\theta)}{\kn_1^2\kn_2^2}=\\[9pt]
     &=&\Frac{\sin\omega%
             \left[-2\sin\gamma\Sam{}\Sbm{}+
                    (\theta{-}\alpha)\Sbm{2}+ (\beta{-}\theta)\Sam{2}
            \right] }%
      {2\sin\gamma\Sam{2}\Sbm{2}} .
   \end{array}
}
We have to prove that the above bracketed expression is
of constant sign. Applying one more substitution, 
\equa{%
   \theta\ieq\omega-x,  \mbox{~~~i.e.~~}
   \alpha-\theta=x+\gamma,\qquad  
   \beta-\theta=x-\gamma,  
}
denote this expression as $F(x;\gamma)$:
\equa{%
   \begin{array}{lcl}
       F(x;\gamma)&=&
           -2\sin\gamma\sin\ffrac{x{+}\gamma}\sin\ffrac{x{-}\gamma}
           -(x{+}\gamma)\sin^2\ffrac{x{-}\gamma}
           +(x{-}\gamma)\sin^2\ffrac{x{+}\gamma}\\[8pt]
                  &=&
            \cos x(\gamma\cos\gamma{+}\sin\gamma) +
            x\sin x\sin\gamma - \sin\gamma\cos\gamma - \gamma .
   \end{array}
}
Since $\theta\In[-\alpha,\beta]\In[-\pi,\pi]$ and $|\omega|\ilt\pi$,
it is sufficient to explore the interval $x\In[-2\pi,2\pi]$.
Because $F(x;\gamma)$ is even with respect to~$x$ 
\ $[F(-x;\gamma)\ieq F(x;\gamma)\,]$,
and odd with respect to the parameter
\ $[F(x;-\gamma)\ieq {-}F(x;\gamma)\,]$,
we explore its behavior only for \GT{\gamma} and \GE{x}.
The plot of $F(x;\frac{2}{3}\pi)$
is shown on the left side of \Reffig{figure14}.
To find extrema of $F(x;\gamma)$ solve the equation 
\EQ{\D{F_x}(x;\gamma)}:
\equa{%
      x\cos x\sin\gamma - \gamma\sin x\cos\gamma = 0.
}
Its non-negative roots \ $x_0,x_1,x_2,\ldots$ \ are: \EQ{x_0}, and
those of the equation $x\cot x\ieq \gamma\cot\gamma$;
in particular, $x_1\ieq\gamma$.
The roots are shown as dots in the right side of \Reffig{figure14},
where the function $x\to x\cot x$ is plotted.
From piecewise monotonicity of this function
it is clear that $x_1\In(0,\pi)$,
$x_2\In(\pi,2\pi)$, etc. We can now describe the
behavior of $F(x;\gamma)$ in $x\In[0,2\pi]$ as follows.
At \EQ{x} function has negative local minimum 
\equa{
   F(0;\gamma)=2\sin^2\ffrac{\gamma}(\sin\gamma{-}\gamma)<0,\quad
%
   F^{\prime\prime}_{xx}(0;\gamma)=\sin\gamma{-}\gamma\cos\gamma>0.
}
It increases to the maximum 
\EQ{F(x_1;\gamma)\ieq F(\gamma;\gamma)},
and then decreases to the subsequent minimum at $x\ieq x_2$.
While increasing from $x_2$ to $x_3\In(2\pi,3\pi)$, function 
passes through the boundary $x\ieq2\pi$ of the interval 
under investigation, still remaining negative at this point:
\LT{F(2\pi;\gamma)\ieq F(0;\gamma)}. It is therefore negative
in $[-2\pi,2\pi]$, except two zeros at $x\ieq{\pm}\gamma$.
The derivative $\D{L}(\theta)$ does not change sign;
and $L(\theta)$ is strictly monotonic. \qed{}
\end{pf}
\Lwfig{t}{0pt}{figure14}%

\Topic{Positional inequalities for short spirals}
%
The following theorem generalizes the earlier results
for ``very short'' spirals (theorem~3 in~\cite{Spiral})
and for convex ones (theorem~5 in~\cite{Theorem5}).
\begin{thm}
\label{LenseTheorem}
Short spiral arc is located within its lense. Except endpoints,
the arc has no common points with lense's boundary.
\end{thm}
\begin{pf}
As the point $P(s)$ moves along the curve,
the circular arcs $APB\ieq\Aarc{-\delta(s)}$, containing $P$,
fill continuously the lense (\Reffig{figure15}).
Because \ $-\delta(0)\ieq \alpha$, and
$\delta(s)$ is strictly monotonic (lemma~\ref{DeltaLemma}),
the curve at the very beginning
deviates immediately from the boundary arc $\Aarc{\alpha}$
to the interior of the lense. 
Near the end point the behavior is similar.\,\qed
\end{pf}
\begin{cor}
A short spiral may cut its chord only once; and this occurs
if and only if the tangent angles $\alpha$ and $\beta$
are nonzero and of the same sign.
\end{cor}
More severe  limitation can be derived if boundary curvatures
are known, as shown in \Reffig{figure16}.
For every inner point of a spiral the unique biarc
$\Biarcab{b}{\alpha}{\beta}$
can be constructed. Thus generated subfamily of biarcs fill
{\em bilense\,}, i.e. the region,
bounded by two biarcs, $AT_1B$ and $AT_2B$.
Arcs $AT_1$ and $T_2B$ belong to boundary circles
of curvatures of the enclosed spiral.

Returning to \Reffigs{figure06}a,b,c, this corresponds to projection
of the point $K\ieq(\kn_1,\kn_2)$
onto the hyperbola \EQ{Q}, yielding
two points, $H_1\ieq(\kn_1,g_2)$ and $H_2\ieq(g_1,\kn_2)$.
They provide in turn two biarcs, marked as~$h_1$ and~$h_2$,
and bilense. We are going to prove that
any short spiral, whose boundary parameters
$(\alpha,\beta,\kn_1,\kn_2)$ belong to the closed region
$KH_1H_2K$, is covered by corresponding bilense.
If point $K$ is being moved backwards and upwards to infinity
(or, in the case of decreasing curvature,
forwards and downwards within 
the lower right branch of hyperbola),
biarcs $h_1$ and $h_2$ approache the boundaries of
the lense~\refeq{BiarcLim}, covering all shorts spirals
with given tangents $\alpha,\beta$
and $-\infty\ile \kn_1\ilt \kn_2\ile \infty$.
It was the subject of theorem~\ref{LenseTheorem}.
\begin{defn}
\label{DefBilense}
\rm 
A {\em bilense} $\Bilense{\alpha,\beta,b_1,b_2}$, 
$0\ile b_1 \ilt b_2 \ile\infty$,
generated by a short non-biarc spiral with end conditions
\equa{%
\alpha,\;\beta,\quad
\mbox{such~that}\quad |\omega|\neq\pi,
\quad\mbox{and}\quad
\kn_1= {-}\sin\alpha - b_1^{-1}\sin\omega,\quad
\kn_2=\sin\beta + b_2\sin\omega,
}
is the region,
bounded by two biarcs $\Biarcab{b_1}{\alpha}{\beta}$ and
$\Biarcab{b_2}{\alpha}{\beta}$, namely,
\equa{%
 \Bilense{\alpha,\beta,b_1,b_2} = \{ (x,y):\;
 (x,y)\in\Biarcab{b}{\alpha}{\beta}\,\}.
}
\end{defn}
\begin{figure}[t]
\Bfig{.44\textwidth}{figure15}{}
\hfill
\Bfig{.44\textwidth}{figure16}{}
\end{figure}
Choosing the parent spiral to be non-biarc,
we force $Q$ to be strictly negative,
and avoid bilense of zero width.
Condition $b_1\ilt b_2$ in Def.~\ref{DefBilense}, 
for both increasing and decreasing curvature,
results from
\equa{%
        Q(\kn_1, \kn_2,\alpha,\beta) =
        \left(1-\Frac{b_2}{b_1}\right)\sin^2\omega < 0.
}
\begin{thm}
\label{BilenseTheorem}
All normalized short spirals with fixed boundary tangents $\alpha,\,\beta$,
such that $|\omega|\neq\pi$, and curvature $\kn(s)$, 
such that $\kn_1\ile \kn(s) \ile \kn_2$,
or        $\kn_1\ige \kn(s) \ige \kn_2$, 
are covered by the corresponding bilense
$\Bilense{\alpha,\beta,b_1,b_2}$,
\equa{
    b_1=\Frac{-\sin\omega}{\kn_1+\sin\alpha},\quad
    b_2=\Frac{\kn_2-\sin\beta}{\sin\omega}\,.
}
\end{thm}
\begin{pf}
\Reffig{figure16} clarifies the proof, based on the monotonicity 
of the map
\equa{%
 \mbox{point on the curve}\;\rightarrow\;
 \mbox{biarc through this point,}
}
i.e. monotonicity of the function $b(s)\ieq b(x(s),y(s))$,
defined by~\refeq{Bxy}. 
Consider the case of increasing curvature.
Denote $C\ieq Z(\bar{s})$ the point
where the spiral meets the bisector of the lense.
From the proof of theorem~\ref{LenseTheorem} it is clear
that such point exists, is unique, and the subarc $AC$ of the spiral
is located in the upper half of the lense.
From $\delta(\bar{s})=-\gamma$ derive equality 
$\omega_1(\bar{s})\ieq\omega_2(\bar{s})$:
\equa{%
 \begin{array}{rcl} 
   0&=&2[\delta(\bar{s}){+}\gamma]
     = 2\mu_2(\bar{s}){-}2\mu_1(\bar{s})+\alpha{-}\beta =\\[6pt]
    &=&[\underbrace{\alpha{-}\mu_1(\bar{s})}_{\alpha_1} +
        \underbrace{\tau(\bar{s}){-}\mu_1(\bar{s})}_{\beta_1}] 
     - [\underbrace{\tau(\bar{s}){-}\mu_2(\bar{s})}_{\alpha_2} +
        \underbrace{\beta{-}\mu_2(\bar{s})}_{\beta_2}]
     = 2[\omega_1(\bar{s})-\omega_2(\bar{s})].
 \end{array}  
}
Denote $\omega_0 = \omega_1(\bar{s})\ieq\omega_2(\bar{s})$ and apply
inequality~\refeq{o1o2}, taking into account that
$0\ilt\omega\ilt\pi$:
\equa{%
   \Omega(\bar{s}) = 2\omega_0{-}\omega < 0
   \So
   \omega_0 < \Frac{\omega}{2} < \Frac{\pi}{2}\,.
}
The map $b(s)\ieq b(x(s),y(s)\,)$ for the points of $AC$
is determined by the first expression of~\refeq{Bxy}.
Applying substitutions
\equa{%
    \begin{array}{l}
        h_1\ieq\sqrt{(x{+}1)^2+y^2},\qquad
        D\ieq(1{-}x^2{-}y^2)\sin\alpha-2y\cos\alpha \ieq \Frac{h_1^2\sin\omega}{b},\\
        \D{x}=\cos\tau, \quad        \D{y}=\sin\tau,    \qquad
        x{+}1 = h_1\cos\mu_1,\quad    y = h_1\sin\mu_1
    \end{array}
}
($x$, $y$, $h_1$, $\mu_1$, $\tau$, $D$, $b$
are functions of $s$),
calculate derivative $\D{b}(s)$:
\equa{%
  \begin{array}{rcl}
    \Dfrac{b}{s} &=&
%
    \Dfrac{}{s}\:\Frac{h_1^2\sin\omega }{D} 
    =2\sin\omega
     \Frac{[(x{+}1)\D{x}{+}y\D{y}]D+
            h_1^2[(x\D{x}{+}y\D{y})\sin\alpha+\D{y}\cos\alpha]}%
          {D^2} =\\[9pt]
     &=&
    2\sin\omega
     \Frac{[(x{+}1)^2-y^2]\sin(\alpha{+}\tau) - 2y(x{+}1)\cos(\alpha{+}\tau)}%
          {D^2} =\\[9pt]
     &=&
    2\sin\omega 
     \Frac{h_1^2 [\cos2\mu_1\sin(\alpha{+}\tau) - \sin2\mu_1\cos(\alpha{+}\tau)]}%
          {[h_1^2\,b^{-1}\sin\omega]^2} =
    \Frac{2 b^2 \sin 2\omega_1]}{h_1^2\,\sin\omega} \ge 0.
  \end{array}
}
This expression is non-negative in $(0,\bar{s}]$
because $\omega_1(s)$ is monotonic increasing with $s$ up to
the value $\omega_1(\bar{s})\ieq \omega_0\ilt\frac{\pi}{2}$.
It may be zero while \EQ{\omega_1(s)}, i.e.
in the case of initial circular subarc,
partially coincident with the initial curvature element of the spiral
(arc $AT_1$).
The value $b(0)$ can be calculated from expansions  
\equa{
  \begin{array}{l}
     x(s)=-1+s\cos\alpha-\Frac{\kn_1}{2}s^2\sin\alpha+O(s^3),\quad
     y(s)=s\sin\alpha+\Frac{\kn_1}{2}s^2\cos\alpha+O(s^3):\\
   b(0)=\lim\limits_{s\to0}b(s)=
  \lim\limits_{s\to0}\Frac{(4+\kn_1^2s^2)\sin\omega}%
                          {-4(\kn_1{+}\sin\alpha)-\kn_1^2s^2\sin\alpha} =
    \Frac{-\sin\omega}{\kn_1+\sin\alpha} = b_1.  
 \end{array}   
}
Similarly, from the second expression of~\refeq{Bxy},
$b(S)\ieq b_2$, and,
due to $2\omega_2(s)\ilt\pi$ in $[\bar{s},S)$, 
the derivative remains non-negative in $[\bar{s},S)$:
\equa{%
   \D{b}(s) =\Dfrac{}{s}\:%
               \Frac{(1{-}x^2{-}y^2)\sin\beta+2y\cos\beta}%
                    {\sin\omega[(x{-}1)^2+y^2]} = \ldots%
    = \Frac{2\,\sin 2\omega_2}{h_2^2\,\sin\omega}\ge 0
}
For decreasing curvature functions $\omega_{1,2}(s)$ as well as the
constant $\omega$ change sign; $b(s)$ remains monotonic increasing. \qed{}
\end{pf}
%

\Reffig{figure17} illustrates an application of these results
to a spiral in the whole. 
A spiral is presented in \Reffig{figure17}a by a set of points 
$P_1,\ldots,P_n$, $n\ieq 11$,  
and tangents $\tau_1,\,\tau_n$ at the endpoints.
The constraint was imposed that the subarcs
$\Arc{P_iP}_{i{+}1}$ were one-to-one projectable
onto the corresponding chord (i.e. $|\acum_i,\bcum_i|\ilt\pi/2$).
For the practice of curves interpolation
it is a quite weak limitation. 
\Lwfig{t}{0pt}{figure17}%

In \Reffig{figure17}b circular arcs $A_1,\ldots,A_n$ are
constructed as follows: 
arc $A_1$ passes through points $P_1$ and $P_2$, 
matching at $P_1$ given tangent $\tau_1$; 
arcs $A_i$, $i\ieq 2,\ldots,n{-}1$ pass
through three consecutive points $P_{i{-}1},P_i,P_{i{+}1}$;
arc $A_n$ passes through two points $P_{n{-}1},P_n$,
matching given tangent $\tau_n$ at the endpoint.
Thus on each chord $P_iP_{i+1}$ we get a lense, bounded by arcs
$A_i$ and $A_{i{+}1}$. The following was proven  
in \cite{Spiral}, theorem~5:
\begin{itemize}
\item[$\bullet$] 
{\em The sequence $k_1,\ldots,k_n$ of curvatures of arcs $A_i$ 
is monotonic.}
\item[$\bullet$]
{\em The union of such lenses covers all spirals,
matching given interpolation data.}
\end{itemize}

The width of this region is of the order $O(h_i^3)$,
\ $h_i\ieq |P_iP_{i{+}1}|$. 
\Reffig{figure17}c shows this construction for
15 points.
The influence of discretization is also seen from comparing
left and right branches of the spiral. 
Thus, the measure of determinancy of a spiral by
inscribed polygonal line is provided, 
without invoking any empirics of particular interpolational
algorithms. 

\Topic
{Existence theorems}
%
The converse of Vogt's theorem, namely, the problem of
joining two line or curvature elements by a spiral arc,
was considered by A.\,Ostrowski~\cite{Ostrowski}.
His solution concerns only $C^2$-continuous convex spirals.
It states that for two given line elements
Vogt's theorem is the sufficient condition for
the existence of such spiral. If curvatures $R_1^{-1}$ and $R_2^{-1}$ 
at the endpoints are involved, additional condition
$D\ilt|R_2{-}R_1|$ is required, $D$~being the distance 
between the centres of two boundary circles of curvature.
Rewritten in terms of this article,
this condition is the particular case of inequality 
\LT{Q(\Kl{1},\Kl{2})}.
Theorem~2 in~\cite{Spiral} establishes this condition
for ``very short'' spirals, regardless of convexity,
and includes the biarc
case as the unique solution if \EQ{Q(\Kl{1},\Kl{2})}.
\begin{thm}
\label{ExistenceShort}
The necessary and sufficient conditions for the existence of
a short spiral curve, matching at the endpoints two given 
curvature elements~\refeq{K1K2c1},
are: modified Vogt's theorem~\refeq{VogtShort} 
and inequality \LE{Q(\Kl{1},\Kl{2})};
if \EQ{Q}, biarc is the unique solution.
\end{thm}
\begin{pf}
Theorems \ref{ShortVogtTheorem} and \ref{MainTheorem} 
prove the necessity of these conditions. 
To prove sufficiency, construct a smooth 
three-arc spiral curve
whose boundary curvature elements are \Kl{1} and \Kl{2}, 
shown as $AC$ and $DB$ in \Reffig{figure18}. 
Apply inversion about the circle 
\equa{%
   \Kl{}^{\star}=\Aarc{\gamma^{\star}}= 
   \Kr{-1,0,\gamma^{\star},\kn^{\star}},\qquad 
   \gamma^{\star}=\gamma/2,\quad
   \quad\kn^{\star}=-\sin(\gamma^{\star}),
}
shown dotted-dashed; its centre is the point 
$O\ieq(0,y^\star)\ieq(0,-\cot\gamma^\star)$.
The inversion is chosen to make the lense symmetric, 
its former bisector
$\Aarc{\gamma}$ is transformed into the chord $\Aarc{0}$. 
For \EQ{\gamma} this is just a symmetry about
X-axis. If $\gamma\neq0$ (which means $|\omega|\ilt\pi$), 
we have to verify that
all the points of the interior of the lense remain
in the interior of its image, i.e. 
$O$ is located outside the lense. 
To do it, we enter lense's boundary ordinates as 
$y_1\ieq\tan\frac{\alpha}{2}$ and 
$y_2\ieq{-\tan\frac{\beta}{2}}$, and check the sign of the product
\equa{
  (y^\star-y_1)(y^\star-y_2)=
\left(-\cot\Frac{\gamma}{2} - \tan\Frac{\alpha}{2}\right) 
 \left(-\cot\Frac{\gamma}{2} + \tan\Frac{\beta}{2}\right) 
 =\Frac{\cos^2\frac{\omega}{2}}%
  {\sin^2\frac{\gamma}{2} \cos\frac{\alpha}{2} \cos\frac{\beta}{2}}>0.  
} 
\Lwfig{t}{0pt}{figure18}%
\noindent New boundary angles
$\D{\alpha},\:\D{\beta}$
can be calculated from conditions
$ \gamma^{\star}\ieq\frac{1}{2}(\alpha{+}\D{\alpha})$,
$-\gamma^{\star}\ieq\frac{1}{2}(\beta{+}\D{\beta})$;
and new curvatures $\D{\kn_1},\:\D{\kn_2}$, from Prop.~\ref{InvCurvatureProp}:
\equa{%
  \begin{array}{lcll}
    \D{\kn_1}&=&           
            2\kn^{\star}(1{-}2Q_{01})-\kn_1,\quad&
            Q_{01}=Q(\Kl{}^{\star},\Kl{1})
                  =\sin^2[(\alpha{-}\gamma^{\star})/2], \\[3pt]
    \D{\kn_2}&=& 
            2\kn^{\star}(1{-}2Q_{02})-\kn_2, &
            Q_{02}=Q(\Kl{}^{\star},\Kl{2})
                  =\sin^2[(\gamma^{\star}{+}\beta)/2]
  \end{array}
}
(to calculate $Q$'s the last equation of~\refeq{Qabc} was used).
This yields
\equa{%
 \D{\alpha}=\D{\beta}=\D{\omega}= -\omega,
 \qquad
 \D{\kn_1}= [-\kn_1{-}\sin\alpha]+\sin\omega,\quad
 \D{\kn_2}=-[\kn_2{-}\sin\beta] - \sin\omega.
}
(e.g., $\D{\kn_1}\ieq {-}2\sin\gamma^\star \cos(\alpha{-}\gamma^\star){-}\kn_1 
 \ieq{-}\sin\alpha{-}\sin(2\gamma^\star{-}\alpha) {-}\kn_1 
 =-\kn_1{-}\sin\alpha{+}\sin\omega$).
For the initial conditions, corresponding to 
increasing curvature (i.e. $0\ilt\omega\ile\pi$, $\kn_1\ilt\kn_2$), 
$\sin\omega$ is non-negative, and the bracketed terms are 
positive. The latter is direct consequence of inequality
\LE{Q(\Kl{1},\Kl{2})} (Cor.~\ref{k1k2Cor}).
New end conditions correspond to decreasing curvature
with boundary curvatures opposite in sign
for whichever signs of~$\kn_{1,2}$: $\D{\kn_1}\igt0\igt\D{\kn_2}$.
Inversed curvature elements are shown as $AC_1$ and $D_1B$.
Show that point $C_1\ieq(x_1,0)$ is located to the left of $D_1\ieq(x_2,0)$:
\equa{%
  \begin{array}{ll}
    AC_1=-\Frac{2\sin\D{\alpha}}{\D{\kn_1}}=
           \Frac{2\sin\omega}{\sin\omega-\kn_1-\sin\alpha},\quad&
    x_1=-1+AC_1=\Frac{\sin\omega+\kn_1+\sin\alpha}{\sin\omega-\kn_1-\sin\alpha},\\
    & \\
    D_1B=\HM{}\Frac{2\sin\D{\beta}}{\D{\kn_2}}=
           \Frac{2\sin\omega}{\sin\omega+\kn_2-\sin\beta},\quad&
    x_2=\HM{}1-D_1B=\Frac{\kn_2-\sin\beta-\sin\omega}{\sin\omega+\kn_2-\sin\beta}.
  \end{array}
}
The denominators in the expressions for $x_{1,2}$ being positive,
it is easy to check that the condition $x_1\ile x_2$ 
is equivalent to \LE{Q}.
The equalities, if occur, are simultaneous,
and two arcs form a unique biarc solution. 
Otherwise the existence of straight line   
$L\ieq\Kr{x_0,0,\lambda,0}$,
smoothly joining two given arcs, is evident; 
its parameters $x_0$ and $\lambda$ 
($x_1\ilt x_0\ilt x_2$, \GT{\lambda})  
can be calculated from two equations \EQ{Q(L,\Kl{1,2})}.
The backward inversion resets the increasing curvature and
yields the sought for solution~--- intermediate arc $L_1$,
image of $L$.
\qed
\end{pf}
\begin{thm}
\label{ExistenceTheorem}
The necessary and sufficient condition for the existence of
a non-biarc spiral curve, matching at the endpoints two given 
curvature elements~\refeq{K1K2cc}, is 
\equa{%
      Q(\Kl{1},\Kl{2})  
      =(k_1c+\sin\alpha)(k_2c-\sin\beta)
        +\sin^2\Frac{\alpha{+}\beta}{2} 
       < 0.
}
\end{thm}
\begin{pf}
The necessity results from theorem~\ref{MainTheorem}.
To prove sufficiency, a sought for spiral can be constructed as 
a three-arc curve.
This  problem was explored in~\cite{InvInv},
and the existence of solutions, all of them being spirals
iff \LT{Q}, was established. 

A simple proof of sufficiency, alternative to that of~\cite{InvInv}, 
can be proposed. Apply inversion,
bringing two given circles into concentric position
(\Reffig{figure19}).
Condition \LT{Q} means that the circles \Kl{1} and \Kl{2} 
do not intersect, and are not tangent.
Therefore such inversion exists. Denote the images of \Kl{1} and \Kl{2} as
\equa{%
  \begin{array}{lll}
    \DKl{1}=\Kr{x_1,y_1,\D{\alpha},\D{k}_1},\quad&
             x_1=a+\sin\D{\alpha}/\D{k}_1, \quad&
             y_1=b-\cos\D{\alpha}/\D{k}_1, \\
    \DKl{2}=\Kr{x_2,y_2,\D{\beta},\,\D{k}_2},&
             x_2=a+\sin\D{\beta}/\D{k}_2, &
             y_2=b-\cos\D{\beta}/\D{k}_2.
  \end{array}
}
The expressions for $x_{1,2}$ and $y_{1,2}$ assure concentricity with
the common centre $(a,b)$. Recalculate~\refeq{DefQ},
which remains invariant under inversion:
\equa{%
      Q(\DKl{1},\DKl{2})=\Frac{-(\D{k}_2-\D{k}_1)^2}{4\D{k}_1\D{k}_2}=
      Q(\Kl{1},\Kl{2})<0.
}
Negative value of the invariant means that two curvatures
$\D{k}_1$ and $\D{k}_2$ are of the same sign, and two circles \DKl{1} and \DKl{2}
are parallel.
All possible intermediate arcs $T_1T_2$ have the same curvature $k_0$,
and the sequence $\D{k}_1,k_0,\D{k_2}$ is monotonic: 
\equa{%
       \Frac{1}{k_0}=\Frac{1}{\D{k_1}}+\Frac{1}{\D{k_2}},\qquad
       k_0=\Frac{2\D{k}_1\D{k}_2}{\D{k}_1+\D{k}_2} \So
       \D{k}_1\lessgtr k_0 \lessgtr \D{k_2}.
}  
An intermediate arc can be constructed for any images 
$A$ and $B$ of two given endpoints.
The backward inversion restores the initial type of monotonicity 
of curvature.\qed{}
\end{pf}
\Lwfig{t}{.8\textwidth}{figure19}%
Note that the strict form of inequality in theorem~\ref{ExistenceTheorem}
is strong enough to exclude both the chord of zero length
and the equality $k_1\ieq k_2$:
with \EQ{c} we get $Q\ieq\sin^2\gamma$, and with $k_1c\ieq k_2c\ieq\kn$ 
the invariant $Q$ takes non-negative form~\refeq{k1eqk2}.
We need not to invoke Vogt's theorem:
if the sum $\alpha{+}\beta$ does not suit
condition $k_1 \lessgtr k_2$~\refeq{VogtAK},
cumulative angles
$\acum\ieq\alpha{\pm}2\pi$ or
$\bcum\ieq\beta{\pm}2\pi$ resolve the contradiction,
resulting to a long spiral as a solution.

\Topic{Conclusions}
%
The author's looking at the theory of spirals
and revisiting Vogt's theorem was initially tied to
the hypothesis that, under certain constraints, the curve  
can be fairly well reproduced from
the interpolation data or inscribed polygon. 
And, contrary to traditional
treatment of this problem, it was interesting to
impose more fundamental constraints than, say, artificial
limits for derivatives, etc. Such constraint was
suggested by the Four-Vertex theorem, and the definition of 
the problem sounded like:

{\em 
Amongst all curves matching given interpolation data
select those having a minimum of vertices; estimate the
region covered by them.
}

The solution for spiral curves, 
i.e. with the minimum of vertices being zero,
was cited here as \Reffig{figure17}.
The preliminary extension for non-spiral curves
was demonstrated in~\cite{Spiral}. Constraints for non-spirals, 
similar to \refeq{k1k2},
were discussed in~\cite{Vertex}.

The tendency to minimize the number of vertices
is very close to the notion of {\em fair curve},
widely discussed in
Computer-Aided Design (CAD) applications~\cite{Fairness}.
 
In recent years much attention is given to curves
with monotonic curvature in CAD related publications.
Refs.~\cite{Theorem5,SpiralSeg} provide just start points
and titles for the bibliographical search.
A lot of research is aimed to extract spiral
subarcs from B\'{e}zier or NURBS curves,
whose polynomial nature is far short of spirality.
We have chosen a somewhat different approach,
having focused on general properties,
induced by monotonicity of curvature.

Most of the earlier studies of spiral curves
seem to have been aimed at obtaining the results similar to the 
Four-Vertex-theorem, and limited to this objective. 
They required continuity of evolute as the basis for the proofs,
and therefore were restricted to the curves with continuous
curvature of constant sign.
Similarly, a lot of CAD applications propose separate
treatment of ``C-shaped'' and ``S-shaped'' spirals,
which is usually unnecessary.

As a final note, we give some attention to the term {\,\em spiral\,}
by itself. Literature treatments are not rigorous
and vary considerably. A curve with a monotonic polar equation
$r(\varphi)$ is often meant by a spiral~(\cite{EncDictMath}, p.~325).
Considering spirality as a property of shape,
the relation to a specific coordinate system is the drawback of 
this definition. Guggenheimer's treatment of spirality
as monotonicity of curvature 
is purely shape related.
Three examples illustrate some ambiguities:
\begin{itemize}
\item[$\bullet$] 
Fermat spiral, $r\ieq a\sqrt{\varphi}$: curvature is not monotonous; 
\item[$\bullet$] 
Cornu spiral $k(s)\ieq s/a^2$: no polar equation for the curve 
as the whole;
\item[$\bullet$] 
C\^otes' spiral $r\ieq a/\cos(k\varphi)$: neither definition is
applicable. 
\end{itemize}
Curves with monotonic curvature comprise an important subset
of planar curves.
The list of  their properties, compiled from~\cite{Vogt},
\cite{Guggen} (pp.~48--54), and this article,
is far short of being complete.
This class of curves is worthy of definite naming.
To designate them, the term {\em ``true spirals''} 
is being proposed.

\paragraph{Acknowledgement.}
The author is grateful to prof. Victor Zalgaller,
who explored his collection of references, 
and called author's attention to the research of \hbox{Wolfgang} Vogt.

\newpage

\newcommand{\Jnum}[1]{$N$#1}
\newcommand{\Jitem}[6]{
{\it#1}, #2. #3, {\bf#4}(#5), #6.}


\begin{thebibliography}{99}
%
\bibitem{InvDist}
\Jitem{Coxeter~H.S.M.}%
{Inversive Distance}%
{Annali di Matematica Pura ed Applicata, $4$}{71}%
{1966}{73--83}
%
\bibitem{EncDictMath}
Encyclopedic Dictionary of Mathematics. Vol.~1.
The Mit Press, Cambridge, 1977.
%
\bibitem{Guggen}
{\it Guggenheimer~H.W.}, Differential geometry.
Dover Publications, New~York, 1977.
%
\bibitem{Tohoku1}
\Jitem{Hirano~K.}{Simple proofs of Vogt's theorem}%
{Tohoku Math. J.}{47}{1940}{126--128}
%
\bibitem{Tohoku2}
\Jitem{Katsuura~S.}{Ein neuer Beweis des Vogtschen Satzes}%
{Tohoku Math. J.}{47}{1940}{94--95}
%
\bibitem{InvInv}
\Jitem{Kurnosenko~A.}%
{Inversive invariant of a pair of circles (in Russian)}%
{Zapiski nauch. sem. POMI}{261}{1999}{167--186}\\
{\tt(http://dbserv.ihep.su/\~{}pubs/prep1999/99-59-e.htm)}.
%
\bibitem{Spiral}
\Jitem{Kurnosenko~A.}{Interpolational properties of planar
spiral curves (in Russian)}%
{Fund. i priklad. matematika
(Fund. and Appl. Math.)}{7}{2001}{\Jnum{2}, 441--463} \\
{\tt(http://dbserv.ihep.su/\~{}pubs/prep1998/98-9-e.htm)}
%
\bibitem{Vertex}
\Jitem{Kurnosenko A.}{Inequalities on the planar curves
in the vicinity of one or two vertices (in Russian)}%
{Zapiski nauch. sem. POMI}{280}{2001}{194--210}
%
%
\bibitem{Biarc92}
\Jitem{Meek~D.S., Walton~D.J.}%
{Approximation of discrete data by G${}^1$ arc splines}
{Computer-Aided Design}{24}{1992}{\Jnum{6}, 301--306}
%
\bibitem{Theorem5}
\Jitem{Meek~D.S., Walton~D.J.}%
{Approximating smooth planar curves by arc splines}%
{J. of Comp. and Appl. Math.}{59}{1995}{221--231}
%
%
\bibitem{Ostrowski}
\Jitem{Ostrowski~A.}{\"{U}ber die Verbindbarkeit von Linien- und
Kr\"{u}mmungselementen d\"{u}rch monoton gekr\"{u}mmte Kurvebogen}%
{Enseignement Math., Ser.~{\rm2}}{2}{1956}{277--292}
%
\bibitem{Vogt}
\Jitem{Vogt~W.}{\"{U}ber monotongekr\"{u}mmte Kurven}%
{J. reine und angew. Math.}{144}{1914}{239--248}
%
\bibitem{SpiralSeg}
\Jitem{Walton~D.J., Meek~D.S.}%
{Planar G$^2$ curve design with spiral segments}%
{Computer-Aided Design}{30}{1998}{\Jnum{7}, 529--538}
%
\bibitem{Fairness}
\Jitem{Yang~X., Wang~G.}%
{Planar point set fairing and fitting by arc splines}%
{Computer-Aided Design}{33}{2001}{35--43}
%
\bibitem{LongBiarcs}
\Jitem{Yong~J.-H., Hu~S.-M., Sun~J.-G.}%
{A note on approximation of discrete data by G${}^1$ arc splines}
{Computer-Aided Design}{31}{1999}{911--915}
%
\end{thebibliography}
\end{document}